\documentclass[reqno]{amsart}

\usepackage{amsfonts}
\usepackage{amssymb}

\newcommand{\ir}{{I \times \R^n}}

\newcommand{\util}{\tilde{u}}

\newcommand{\lqlr}{{L^q_t L^r_x}}
\def\be#1{\begin{equation} \label{#1}}
\def\bi{\begin{itemize}}
\def\bs{\begin{split}}
\def\es{\end{split}}
\def\ba{\begin{align}}
\def\bas{\begin{align*}}
\def\ea{\end{align}}
\def\eas{\end{align*}}
\def\R{{{\mathbb R}}}
\def\C{{{\mathbb C}}}
\def\Z{{{\mathbb Z}}}

\def\eps{\varepsilon}

\def\emph#1{{\it #1}}
\def\textbf#1{{\bf #1}}

\theoremstyle{plain}
  \newtheorem{theorem}[subsection]{Theorem}
  
  \newtheorem{proposition}[subsection]{Proposition}
  \newtheorem{lemma}[subsection]{Lemma}
  \newtheorem{corollary}[subsection]{Corollary}

\theoremstyle{remark}
  \newtheorem{remark}[subsection]{Remark}

\theoremstyle{definition}

\numberwithin{equation}{section}

\include{psfig}

\begin{document}

\title[Stability of energy-critical NLS]{Stability of energy-critical nonlinear Schr\"odinger equations in high dimensions}
\author{Terence Tao}
\address{Department of Mathematics, UCLA, Los Angeles CA 90095-1555}
\email{ tao@math.ucla.edu}

\author{Monica Visan}
\address{Department of Mathematics, UCLA, Los Angeles CA 90095-1555}
\email{ mvisan@math.ucla.edu} 

\vspace{-0.3in}
\begin{abstract}
We develop the existence, uniqueness, continuity, stability, and
scattering theory for energy-critical nonlinear Schr\"odinger
equations in dimensions $n \geq 3$, for solutions which have large,
but finite, energy and large, but finite, Strichartz norms.  For
dimensions $n \leq 6$, this theory is a standard extension of the
small data well-posedness theory based on iteration in Strichartz
spaces. However, in dimensions $n > 6$ there is an obstruction to this
approach because of the subquadratic nature of the nonlinearity
(which makes the derivative of the nonlinearity non-Lipschitz). We
resolve this by iterating in exotic Strichartz spaces instead. The
theory developed here will be applied in a subsequent paper of the
second author, \cite{monica-thesis}, to establish global
well-posedness and scattering for the defocusing energy-critical
equation for large energy data.
\end{abstract}


\maketitle

\section{Introduction}

We study the initial value problem for the following nonlinear Schr\"odinger equation in spacetime $\R\times\R^n = \R^{}_t \times \R^n_x$, $n\geq 3$,
\begin{equation}\label{equation 1}
\left\{\begin{aligned}
i u_t +\Delta u &= f(u)\\
u(t_0,x) &= u_0(x) \in \dot H^1(\R^n)
\end{aligned}\right.
\end{equation}
where $u(t,x)$ is a complex-valued function in spacetime $\ir$ for some time interval $I$ containing $t_0$, and
the nonlinearity\footnote{The analysis here also extends to systems where $u$ takes values in some finite-dimensional
complex vector space $\C^m$ and $f: \C^m \to \C^m$ obeys analogous estimates to those presented here.  However, we have
elected to only present the scalar case to simplify the exposition.} $f:\C\to \C$ is continuously differentiable
and obeys the power-type estimates
\begin{align}
f(u)&=O\bigl(|u|^{\frac{n+2}{n-2}}\bigr), \label{f}\\
f_z(u), \ f_{\bar{z}}(u) &= O\bigl(|u|^{\frac{4}{n-2}}\bigr),\label{f_z}\\
\left.\begin{aligned}
f_z(u)-f_z(v) \\
f_{\bar{z}}(u)-f_{\bar{z}}(v)
\end{aligned}\right\} &= \label{f_z diff}
\begin{cases}
O\bigl(|u-v|^{\frac{4}{n-2}}\bigr), \ \text{if} \ n>6,\\
O\bigl(|u-v|(|u|^{\frac{6-n}{n-2}}+|v|^{\frac{6-n}{n-2}})\bigr), \ \text{if} \ 3\leq n\leq6,
\end{cases}
\end{align}
where $f_z$, $f_{\bar{z}}$ are the usual complex derivatives
$$f_z := \frac{1}{2}\Bigl(\frac{\partial f}{\partial x} - i\frac{\partial f}{\partial y}\Bigr), \quad f_{\bar{z}} :=
\frac{1}{2}\Bigl(\frac{\partial f}{\partial x} + i\frac{\partial f}{\partial y}\Bigr).$$
For future reference, we observe the chain rule
\begin{equation}\label{chain}
\nabla f(u(x)) = f_z(u(x)) u(x) + f_{\bar z}(u(x)) \overline{u(x)}
\end{equation}
as well as the closely related integral identity
\begin{equation}\label{ftoc}
f(u)-f(v)
= \int_0^1 \bigl[f_{z}\bigl(v+\theta(u-v)\bigr)(u-v)+f_{\bar{z}}\bigl(v+\theta(u-v)\bigr)\overline{(u-v)}\bigr]d\theta
\end{equation}
for any $u,v \in \C$; in particular, from \eqref{f_z}, \eqref{ftoc},  and the triangle inequality, we have the estimate
\begin{align}\label{f diff}
\bigl|f(u)-f(v)\bigr|\lesssim |u-v|\bigl(|u|^{\frac{4}{n-2}}+|v|^{\frac{4}{n-2}}\bigr).
\end{align}
Following \cite{cazenave:book}, we shall only consider \emph{strong solutions}\footnote{In dimensions $n>4$,
the nonlinearity is unlikely to be smooth and so we cannot use the usual trick of working entirely with classical
(Schwartz class) solutions first and then taking limits, unless one also performs some regularization of the nonlinearity.
If however the nonlinearity is smooth, then it is easy to see that Schwartz initial data leads to Schwartz solutions
and, by using the well-posedness and stability theory which we will develop, one can then express the strong solutions
given here as the unique strong limit of classical solutions.} to \eqref{equation 1}, by which we mean
solutions $u \in C^0_t \dot H^1_x(\ir)$ to the integral (or Duhamel) formulation
\begin{equation}\label{iteration}
u(t) = e^{i(t-t_0)\Delta} u_0 - i \int_{t_0}^t e^{i(t-s)\Delta} f(u(s))\ ds
\end{equation}
of the equation (we will explain our notation more fully in the next section). Note that
by combining Sobolev embedding, $\dot H^1(\R^n) \subseteq L^{2n/(n-2)}(\R^n)$, and
\eqref{f} we see that $f(u(s))$ will be a tempered distribution uniformly in $s$.  Consequently,
there is no difficulty interpreting \eqref{iteration} in the sense of tempered
distributions at least.

The most important examples of nonlinearities of the above type are the defocusing power
nonlinearity $f(u) := + |u|^{\frac{4}{n-2}} u$ and the focusing power nonlinearity
$f(u) := -|u|^{\frac{4}{n-2}} u$.  Note that in these cases, the problem \eqref{equation 1} is invariant under the scaling
$$ u(t,x) \mapsto \frac{1}{\lambda^{(n-2)/2}} u\Bigl(\frac{t}{\lambda^2}, \frac{x}{\lambda}\Bigr),
 \quad u_0(x) \mapsto \frac{1}{\lambda^{(n-2)/2}} u_0\Bigl(\frac{x}{\lambda}\Bigr)$$
and that this scaling also preserves the $\dot H^1(\R^n)$ norm of $u_0$.  Thus we shall view \eqref{equation 1} as an
\emph{$\dot H^1_x$-critical} or \emph{energy-critical} equation\footnote{In this paper we shall use \emph{energy}
synonymously with (the square of) the $\dot H^1_x$ norm; the potential energy (which is essentially the $L^{2n/(n-2)}_x$ norm) will not play a role in our discussions.}.

The Cauchy problem \eqref{equation 1}
has been studied extensively in the literature, especially in dimensions $n=3,4$; see the references below
and particularly the books \cite{borg:book}, \cite{cazbook}, \cite{cazenave:book}.  One can divide the theory
into two parts: the ``local'' theory in which the solution is either restricted to have small energy, or to have
a certain spacetime norm small (which can be achieved for instance by localizing the time interval) and
the ``global'' theory in which there is no size restriction on the solution (other than finite energy).
Our focus here will be on the local theory, in which the exact choice of the nonlinearity $f$ does not play a major role;
in particular, there will be no distinction between the focusing and defocusing cases.  However, the results here will
be applied towards the global theory. Specifically, they will be used (together with several other tools,
notably an interaction Morawetz inequality and Bourgain's induction on energy argument) in a future paper of the
second author, \cite{monica-thesis}, establishing global well-posedness and scattering for the defocusing energy-critical
equation for large $\dot H^1(\R^n)$ data for dimensions $n > 4$; the cases $n=3,4$ were treated in \cite{gopher},
\cite{rv} respectively (see also the  works \cite{borg:scatter}, \cite{borg:book}, \cite{grillakis:scatter},
\cite{tao:gwp radial} treating the radial case).

For the energy-critical local theory it is convenient to introduce a number of scale-invariant
function spaces.  We use $L^r_x(\R^n)$ to denote the Banach space of functions $f: \R^n \to \C$ whose norm
$$ \| f \|_r := \Bigl(\int_{\R^n} |f(x)|^r\ dx\Bigr)^{1/r}$$
is finite (with the usual modification when $r=\infty$). Similarly, for any spacetime slab $\ir$, we use $L^q_t L^r_x(\ir)$ to denote the
Banach space of functions $u: \R^n \to \C$ whose norm
$$ \| u \|_{L^q_t L^r_x(\ir)} = \|u\|_{q,r} := \Bigl(\int_I \|u(t)\|_r^q\ dt\Bigr)^{1/q}$$
is finite (with the usual modification when $q=\infty$).  We will omit mention of the slab $\ir$ when it is clear from context.
In addition to the energy space $C^0_t \dot H^1_x(\ir)$, we will also need the
Strichartz space $\dot W = \dot W(\ir)$, defined on $\ir$ as the closure of the test functions under the norm
\begin{equation}\label{normbounded}
 \|u\|_{\dot W} := \|\nabla u\|_{\frac{2(n+2)}{n-2}, \frac{2n(n+2)}{n^2+4}}
\end{equation}
which is a convenient norm that is particularly well adapted for controlling solutions to \eqref{equation 1}.  Similarly,
we introduce the nonlinearity space $\dot N^1(\ir)$ defined on $\ir$ as the closure of the test functions under the norm
\begin{equation}\label{n1-def}
 \| F \|_{\dot N^1} := \| \nabla F \|_{2, \frac{2n}{n+2}}.
\end{equation}
This space is useful for controlling the forcing term $f(u)$.

A large part of the local theory for these energy-critical equations was worked out by
Cazenave and Weissler, \cite{cw0}, \cite{cwI}, building upon earlier work of Ginibre-Velo,
\cite{gv:localreference}, and Kato, \cite{kato}, for the energy-subcritical case.  In this
work, global solutions were constructed for small energy data and local solutions were
constructed for large energy data, though, as is to be expected for a critical equation,
the time of existence depends on the profile of the initial data and not simply on the
energy.  Furthermore, these solutions $u$ were unique in a certain Strichartz space and
in this space the solution depended continuously\footnote{Strictly speaking, the
continuity in $C^0_t H^1_x$ was only demonstrated for the defocusing equation; for the
general equation the continuity was established in $L^q_t H^1_x$ for any finite $q$.  See
\cite{cazenave:book} for further discussion.  In this paper we shall obtain the stronger
result of H\"older continuity for the general equation in the scale-invariant space
$C^0_t \dot H^1_x$.} on the initial data in the energy space $\dot H^1(\R^n)$.
A later argument of Cazenave, \cite{cazenave:book}, also demonstrates that the uniqueness
is in fact unconditional in the category of strong solutions (see also \cite{katounique},
\cite{twounique}, \cite{FPT_NLSunique}, \cite{gopher} for some related arguments; we
reproduce the argument in Proposition \ref{unconditional} below).

We now give an extremely oversimplified sketch of how these results are obtained.  We
rely heavily on the spaces $\dot W(\ir)$ and $\dot N^1(\ir)$ introduced earlier. The
fundamental tools are the Strichartz estimate\footnote{Here and in the sequel, $X \lesssim Y$
will denote an estimate of the form $X \leq C Y$ where $C > 0$ depends only on
the dimension $n$.}
\begin{equation}\label{strich-basic}
\bigl\| e^{i(t-t_0)\Delta} u_0 + \int_{t_0}^t e^{i(t-s)\Delta} F(s)\ ds \bigr\|_{\dot W}
\lesssim \| u_0 \|_{\dot H^1(\R^n)} + \| F \|_{\dot N^1},
\end{equation}
whenever the right-hand side is finite (and on any spacetime slab $\ir$ containing $t=t_0$), as well as the Sobolev embedding
\begin{equation}\label{sob-embed}
 \| u \|_{\frac{2(n+2)}{n-2},\frac{2(n+2)}{n-2}} \lesssim \| u \|_{\dot W}
 \end{equation}
for all $u \in \dot W(\ir)$.  We also take
advantage of the H\"older inequality
\begin{equation}\label{holder}
\| v^{4/(n-2)} \nabla u \|_{2,\frac{2n}{n+2}}
    \leq \| v \|_{\frac{2(n+2)}{n-2},\frac{2(n+2)}{n-2}}^{4/(n-2)} \| u \|_{\dot W},
\end{equation}
which in conjunction with \eqref{chain} and \eqref{sob-embed} implies that
\begin{equation}\label{chain-game}
\| f(u) \|_{\dot N^1} \lesssim \| u \|_{\frac{2(n+2)}{n-2},\frac{2(n+2)}{n-2}}^{4/(n-2)} \| u \|_{\dot W} \lesssim \|u\|_{\dot W}^{(n+2)/(n-2)}.
\end{equation}
Combining this with \eqref{strich-basic} we see that
strong solutions to \eqref{iteration} which lie in $\dot W(\ir)$ obey the \emph{a priori} estimate
$$
\| u \|_{\dot W} \lesssim \| u_0 \|_{\dot H^1_x} + \| u \|_{\dot W}^{(n+2)/(n-2)},
$$
which then suggests that $u$ stays small in $\dot W(\ir)$ whenever $u_0$ is small in $\dot H^1(\R^n)$.
This can be made more rigorous by setting up an iteration scheme to construct $u$; the case of
large energy can be dealt with by an appropriate truncation of time (to reduce the $\dot W(\ir)$ norm of
$e^{i(t-t_0)\Delta} u_0$).  The continuous dependence on the data is not difficult in rough spaces (e.g.
in $C^0_t L^2_x(\ir)$), but to obtain continuous dependence in an energy-critical space was previously
only obtained in the defocusing case by a delicate argument, requiring the energy conservation law and a
sharp form of Fatou's lemma; see \cite{cazenave:book} for details.

The above results (which were obtained by iteration in Strichartz spaces) already form a quite
satisfactory local well-posedness theory in the energy space for the above equations.  However, there
are still some points that need to be resolved. Firstly, the known arguments that establish
continuous dependence on the data do not necessarily establish \emph{uniformly} continuous dependence on
the data in energy-critical spaces (though they do apply in more supercritical spaces, such as spaces
that scale like $L^2_x(\R^n)$).  In fact, in the low dimensional cases, $n \leq 6$, it is possible to
establish \emph{Lipschitz} dependence on the data, because then we can differentiate equation
\eqref{iteration} using \eqref{chain} to obtain
\begin{equation}\label{diff-iteration}
\nabla u(t) = e^{i(t-t_0)\Delta} \nabla u_0 -i  \int_{t_0}^t e^{i(t-s)\Delta}[ f_z(u(s)) \nabla u(s)
    + f_{\bar{z}}(u(s)) \nabla\overline{ u(s)} ]\ ds
\end{equation}
and rely on the Lipschitz bounds \eqref{f_z diff} to control the difference between two solutions with slightly different
data.  The $n=3,4$ instances of this argument can be found in \cite{gopher}, \cite{rv}, as well as
Section~\ref{short-sec} below.  However, in the high-dimensional case the derivative maps $f_z$, $f_{\bar{z}}$ are
merely H\"older continuous instead of Lipschitz.  If one then tries to analyze differences of two solutions to
\eqref{diff-iteration} with slightly different initial data, one no longer obtains any useful bounds\footnote{This is
basically because any estimate of the form $A \leq \eps + \delta A^\theta$, where $0 < \theta < 1$, does not imply a
bound on $A$ which goes to zero as $\eps$ goes to zero, in contrast to the Lipschitz case, $\theta = 1$, in which one can
obtain a bound of the form $A \leq 2\eps$ (say) if a continuity argument is available and $\delta \leq \frac{1}{2}$.}.

Closely related to the continuous dependence of the data is the \emph{stability} theory for the equation \eqref{equation 1}.  By this we mean
the following type of property: given an \emph{approximate} solution
\begin{equation}\label{equation 1-approx}
\begin{cases}
i \tilde u_t +\Delta \tilde u &= f(\tilde u) + e\\
\tilde u(t_0,x) &= \tilde u_0(x) \in \dot H^1(\R^n)
\end{cases}
\end{equation}
to \eqref{equation 1}, with $e$ small in a suitable space and $\tilde u_0 - u_0$ small in $\dot H^1_x$, is it possible to show that the
\emph{genuine} solution $u$ to \eqref{equation 1} stays very close to $\tilde u$ (for instance, in the $\dot W(\ir)$-norm)?  Note that the question of
continuous dependence of the data corresponds to the case $e=0$.  Again, if $n \leq 6$, an analysis based on subtracting \eqref{diff-iteration}
from the analogous equation for $\tilde u$ and applying Strichartz estimates will yield a satisfactory theory, at least when $\tilde u$
has small $\dot W(\ir)$-norm; the case of large $\dot{W}$-norm can be obtained by partitioning the time interval and iterating the small norm theory.
See \cite{gopher}, \cite{rv} for instances of this argument (which also appears implicitly in \cite{borg:scatter},
\cite{borg:book}).  This type of approach does not work directly in dimensions $n > 6$ as the iteration is only H\"older
continuous instead of Lipschitz and so, one is unable to close the argument effectively, even when one localizes time
to make various norms small.

The purpose of this paper is to complete the previous analysis of Cazenave and Weissler for the Cauchy
problem \eqref{equation 1}, by establishing a local well-posedness and stability theory which is
H\"older continuous in energy-critical spaces and that applies even for large energy data, provided that
the $\dot W(\ir)$-norm is known to be bounded.  This type of result is necessary for induction on energy
type arguments, and will be applied in a subsequent paper of the second author, \cite{monica-thesis}. The
main new tools will be an exotic Strichartz estimate, together with an estimate of fractional chain rule
type; the point of working in an exotic Strichartz space is that it becomes possible to work with a
small fractional derivative rather than a full derivative while still remaining energy-critical with
respect to scaling.  A very similar technique was employed by Nakanishi \cite{nakanishi} for the energy-critical non-linear Klein-Gordon
equation in high dimensions.

We now present our main results.  We begin with a preliminary (and standard) local well-posedness
theorem, which gives existence and uniqueness, as well as Lipschitz continuity, but in rough
(supercritical) topologies. It does not require any exotic Strichartz spaces or nonlinear estimates,
relying instead on iteration in the usual Strichartz spaces, the Leibnitz rule \eqref{chain}, and
H\"older's inequality. It is convenient to place the initial data in the inhomogeneous Sobolev space
$H^1(\R^n)$ rather than the homogeneous one $\dot H^1(\R^n)$; once we obtain the H\"older continuity
estimates on the solution map, we will see that we can easily pass to the homogeneous space from the
inhomogeneous one by a limiting argument.

\begin{theorem}[Standard local well-posedness, \cite{cw0}, \cite{cwI}]\label{lwp}
Let $I$ be a compact time interval that contains $t_0$.
Let $u_0\in H^1(\R^n)$ be such that
\begin{align}
\|e^{i(t-t_0)\Delta}u_0\|_{\dot W(\ir)}\leq \eta \label{small}
\end{align}
for some $0<\eta\leq \eta_0$ where $\eta_0 > 0$ is a small constant. Then there exists a unique solution $u$ to
\eqref{equation 1} in $\dot{S}^1(\ir)$ (this space will be defined in the next section). Moreover, we have the bounds
\begin{align}
\|u\|_{\dot W(\ir)}&\lesssim \eta\label{small solution in L^p}\\
\|u\|_{\dot S^1(\ir)}&\lesssim \| u_0\|_{\dot H^1(\R^n)} + \eta^{\frac{n+2}{n-2}}\label{bounded-s1dot}\\
\|u\|_{\dot S^0(\ir)}&\lesssim \| u_0\|_{L^2(\R^n)}. \label{bounded-s0}
\end{align}
We can generalize \eqref{bounded-s0} as follows:
if $\tilde u_0 \in H^1(\R^n)$ is another initial data obeying the hypothesis \eqref{small}, with associated
solution $\tilde u \in \dot S^1(\ir)$, then we have the Lipschitz bound
\begin{equation}\label{lipschitz}
\| \tilde u - u \|_{\dot S^0(\ir)} \lesssim \| \tilde u_0 - u_0 \|_{L^2(\R^n)}.
\end{equation}
\end{theorem}

For the convenience of the reader we give a proof of this standard result in Section~\ref{lwp-sec}; it
does not require the H\"older continuity hypothesis \eqref{f_z diff}. Remarkably, there is no
restriction on the $H^1(\R^n)$-norm of the initial data, though we do require that this norm is finite.
Instead, we have the smallness condition \eqref{small}.  Note from the Strichartz estimate
\eqref{strich-basic} that
$$
\|e^{i(t-t_0)\Delta}u_0\|_{\dot W(\ir)} \lesssim \| u_0 \|_{\dot H^1(\R^n)} < \infty.
$$
This already gives local existence for all large energy data since, from the above Strichartz estimate and monotone
convergence, we can establish the hypothesis \eqref{small} for intervals $I$ that are sufficiently small (but note
that the size of $I$ will depend on the profile of the initial data $u_0$ and not just on its $\dot H^1(\R^n)$-norm).
Another immediate corollary of the above theorem and the Strichartz estimate is

\begin{corollary}[Global well-posedness for small $H^1_x$ data, \cite{cw0}, \cite{cwI}]\label{cor lwp}
Let $u_0\in H^1_x$ be such that
\begin{align}\label{initial data small}
\|u_0\|_{\dot{H}^1(\R^n)}\leq \eta_0
\end{align}
for some small absolute constant $\eta_0 > 0$ depending only on the dimension $n$.
Then, there exists a unique global $\dot{S}^1(\R \times \R^n)$ solution $u$ to \eqref{equation 1}. Moreover,
\begin{align*}
\|u\|_{\dot W(\R\times\R^n)}, \|u\|_{\dot S^1(\R\times\R^n)} &\lesssim \|u_0\|_{\dot{H}^1(\R^n)} \\
\|u\|_{\dot S^0(\R\times\R^n)} &\lesssim \|u_0\|_{L^2(\R^n)}.
\end{align*}
\end{corollary}

By combining the above standard theory with the exotic Strichartz estimate and fractional chain rule, we obtain our first main result.

\begin{theorem}[Short-time perturbations]\label{short-time theorem}
Let $I$ be a compact time interval and let $\util$ be an approximate solution to \eqref{equation 1} on $\ir$
in the sense that
$$
(i\partial_t+\Delta)\util=f(\util)+e
$$
for some function $e$. Suppose that we also have the energy bound
\begin{align}
\|\util\|_{L_t^{\infty}\dot{H}^1_x(\ir)}&\leq E \label{finite energy-0}
\end{align}
for some constant $E>0$. Let $t_0\in I$ and let $u(t_0)\in H^1(\R^n)$ be close to $\util(t_0)$ in the sense that
\begin{align}\label{close-0}
\|u(t_0)-\util(t_0)\|_{\dot{H}^1(\R^n)}\leq E'
\end{align}
for some $E'>0$. Moreover, assume the smallness conditions
\begin{align}
\|\util\|_{\dot W(\ir)} &\leq \eps_0 \label{finite S norm-0}\\
\Bigl(\sum_N \|P_N e^{i(t-t_0)\Delta}\bigl(u(t_0)-\util(t_0)\bigr)\|^2_{\dot W(\ir)}\Bigr)^{1/2} &\leq \eps \label{closer-0} \\
\|e\|_{\dot N^1(\ir)} &\leq \eps \label{error small-0}
\end{align}
for some $0<\eps\leq \eps_0$, where $\eps_0=\eps_0(E,E') > 0$ is a small constant. Then there exists a solution
$u \in S^1(\ir)$ to \eqref{equation 1} on $\ir$ with the specified initial data $u(t_0)$ at time $t=t_0$ satisfying
\begin{align}
\|u-\util\|_{\dot W(\ir)} &\lesssim \eps+\eps^{\frac{7}{(n-2)^2}} \label{close in L^p-0}\\
\|u-\util\|_{\dot{S}^1(\ir)} &\lesssim E'+\eps+ \eps^{\frac{7}{(n-2)^2}}\label{close in S^1-0}\\
\|u\|_{\dot{S}^1(\ir)} &\lesssim E+E' \label{u in S^1-0}\\
\bigl\|(i \partial_t + \Delta)(u-\util) + e\bigr\|_{\dot N^1(\ir)} &\lesssim \eps+\eps^{\frac{7}{(n-2)^2}}. \label{one}
\end{align}
\end{theorem}

We prove this theorem in Section \ref{short-sec}.
This theorem allows the energy of $\util$ and $u$ to be large as long as the error $e$ is small in a suitable norm, the free evolution of
$u - \util$ is small in another norm, and $\util$ itself is small in a third norm.  The $n=3,4$ cases of this theorem are in \cite{gopher},
\cite{rv} respectively, and the argument there extends easily to dimensions $n=5,6$. However, the cases $n>6$ require
a more delicate argument.
It is probably possible to replace the Besov-norm type expression on the left-hand side of \eqref{closer-0} by the
Sobolev-norm expression $\| e^{i(t-t_0)\Delta} (u(t_0)-\util(t_0)) \|_{\dot W(\ir)}$, but this would presumably require
the Coifman-Meyer theory of paraproducts and we will not pursue it here to simplify the exposition.  The H\"older exponent
$\frac{7}{(n-2)^2}$ can be improved somewhat (perhaps to $\frac{4}{n-2}$) but we will not seek the optimal exponent here
(for applications, all that is important is that this exponent is positive).  However, it seems doubtful that one
can obtain Lipschitz type bounds when the dimension $n$ is sufficiently large.  Indeed, if one had Lipschitz continuous
dependence on the initial data then, by taking variations of \eqref{equation 1} in $u$, one must (formally at least)
have that the linearized equation
$$
\begin{cases}
i v_t +\Delta v &= f_z(u(s)) v(s) + f_{\bar{z}}(u(s)) \overline{v(s)}\\
v(t_0,x) &= v_0(x) \in \dot H^1(\R^n),
\end{cases}
$$
which is a system of linear Schr\"odinger equations in $v$ and $\overline{v}$ with time-dependent, non-self-adjoint potential,
keeps the $H^1_x$ norm of $v(t)$ bounded in time.
If $n$ is large enough, it is unlikely that $f_z(u)$ (which behaves like $|u|^{4/(n-2)}$) remains in $H^1_x$ or even in
$H^{1/2}_x$; thus, it seems that solutions to this equation may leave $H^1_x$ even when local smoothing effects (which typically gain half a derivative of
regularity at most) are exploited\footnote{However, it is still possible to iterate the above linearized equation in lower
regularity spaces of the same scaling as $\dot H^1_x$ and, in particular, using the $X$ and $Y$ spaces introduced in
Section \ref{short-sec}, recover estimates of Lipschitz type.  Indeed, it was this observation for the linearized equation
which eventually led to the arguments here.}.  It seems however to be somewhat difficult to convert these heuristic
arguments into an explicit and rigorous demonstration of non-Lipschitz continuity and we will not pursue this issue here.

By an easy iteration argument (as in \cite{gopher}, \cite{rv}) based on partitioning the
time interval, we can relax hypothesis \eqref{finite S norm-0} by allowing $\util$ to be large (but still bounded in some norm):

\begin{theorem}[Long-time perturbations]\label{long-time theorem}
Let $I$ be a compact time interval and let $\util$ be an approximate solution to \eqref{equation 1} on $I\times\R^n$
in the sense that
$$
(i\partial_t+\Delta)\util=f(\util)+e
$$
for some function $e$. Assume that
\begin{align}
\|\util\|_{L_{t,x}^{\frac{2(n+2)}{n-2}}(\ir)}&\leq M \label{finite S norm} \\
\|\util\|_{L_t^{\infty}\dot{H}^1_x(\ir)}&\leq E \label{finite energy}
\end{align}
for some constants $M, E>0$. Let $t_0\in I$ and let $u(t_0)$ close to $\util(t_0)$ in the sense that
\begin{align}\label{close}
\|u(t_0)-\util(t_0)\|_{\dot{H}^1_x}\leq E'
\end{align}
for some $E'>0$. Assume also the smallness conditions
\begin{align}
\Bigl(\sum_N \|P_N e^{i(t-t_0)\Delta}\bigl(u(t_0)-\util(t_0)\bigr)\|^2_{\dot W(\ir)}\Bigr)^{1/2} &\leq \eps \label{closer} \\
\|e\|_{\dot N^1(\ir)}&\leq \eps \label{error small}
\end{align}
for some $0<\eps \leq \eps_1$, where $\eps_1=\eps_1(E, E', M)$ is a small constant. Then there exists a solution
$u$ to \eqref{equation 1} on $\ir$ with the specified initial data $u(t_0)$ at time $t=t_0$ satisfying
\begin{align}
\|u-\util\|_{\dot W(\ir)}&\leq C(E ,E', M)\bigl(\eps+\eps^{\frac{7}{(n-2)^2}}\bigr) \label{close in L^p}\\
\|u-\util\|_{\dot{S}^1(\ir)}&\leq C(E ,E', M)\bigl(E'+\eps+ \eps^{\frac{7}{(n-2)^2}}\bigr) \label{close in S^1}\\
\|u\|_{\dot{S}^1(\ir)}&\leq C(E, E', M). \label{u in S^1}
\end{align}
Here, $C(E,E',M) > 0$ is a non-decreasing function of $E,E',M$, and the dimension~$n$.
\end{theorem}

We prove this theorem in Section \ref{large-sec}.  As a corollary of this theorem we also
obtain a satisfactory scattering theory for these equations provided that one assumes a
global $L^{2(n+2)/(n-2)}_{t,x}$ bound on solutions; see Corollary \ref{L^p implies
scattering}. This global bound is not difficult to obtain for small energy data (it
follows directly from Theorem \ref{lwp}), but for large energy data the situation is
substantially more delicate and requires further structural information on the
nonlinearity.  In \cite{monica-thesis} the second author will establish this bound for
the defocusing equation; in the focusing case, Glassey's virial identity can be used to
establish blowup for certain large data (even if the data is smooth, spherically
symmetric, and compactly supported); see \cite{glassey}.

\textit{Acknowledgements.} The first author is partly supported by a grant from the Packard Foundation.  The authors also thank Kenji Nakanishi for
pointing out the connections between this paper and \cite{nakanishi}.

\subsection{Notation}

We will often use the notation $X \lesssim Y$ whenever there exists some constant $C$ so that $X \leq CY$. Similarly,
we will write $X \sim Y$ if $X \lesssim Y \lesssim X$.  We say $X \ll Y$ if $X \leq cY$ for some small constant $c$.
The derivative operator $\nabla$ refers to the spatial variable only.

We define the Fourier transform on $\R^n$ to be
$$
\hat f(\xi) := \int_{\R^n} e^{-2 \pi i x \cdot \xi} f(x) dx.
$$

We shall use of the fractional differentiation operators $|\nabla|^s$ defined by
$$
\widehat{|\nabla|^sf}(\xi) := |\xi|^s \hat f (\xi).
$$
These define the homogeneous Sobolev norms
$$
\|f\|_{\dot H^s_x} := \| |\nabla|^s f \|_{L_x^2}.
$$

Let $e^{it\Delta}$ be the free Schr\"odinger propagator.  In physical space this is given by the formula
$$
e^{it\Delta}f(x) = \frac{1}{(4 \pi i t)^{n/2}} \int_{\R^n} e^{i|x-y|^2/4t} f(y) dy,
$$
while in frequency space one can write this as
\begin{equation}\label{fourier rep}
\widehat{e^{it\Delta}f}(\xi) = e^{-4 \pi^2 i t |\xi|^2}\hat f(\xi).
\end{equation}
In particular, the propagator preserves the above Sobolev norms and obeys the \emph{dispersive inequality}
\begin{equation}\label{dispersive ineq}
\|e^{it\Delta}f(t)\|_{L_x^\infty} \lesssim |t|^{-\frac{n}{2}}\|f(t)\|_{L_x^1}
\end{equation}
for all times $t\neq 0$.  We also recall \emph{Duhamel's formula}
\begin{align}\label{duhamel}
u(t) = e^{i(t-t_0)\Delta}u(t_0) - i \int_{t_0}^t e^{i(t-s)\Delta}(iu_t + \Delta u)(s) ds.
\end{align}

We will also need some Littlewood-Paley theory.  Specifically, let $\varphi(\xi)$ be a smooth bump supported in the ball
$|\xi| \leq 2$ and equalling one on the ball $|\xi| \leq 1$.  For each dyadic number $N \in 2^\Z$ we define the
Littlewood-Paley operators
\begin{align*}
\widehat{P_{\leq N}f}(\xi) &:=  \varphi(\xi/N)\hat f (\xi),\\
\widehat{P_{> N}f}(\xi) &:=  (1-\varphi(\xi/N))\hat f (\xi),\\
\widehat{P_N f}(\xi) &:=  [\varphi(\xi/N) - \varphi (2 \xi /N)] \hat f (\xi).
\end{align*}
Similarly, we can define $P_{<N}$, $P_{\geq N}$, and $P_{M < \cdot \leq N} := P_{\leq N} - P_{\leq M}$, whenever $M$ and
$N$ are dyadic numbers.  Sometimes we may write $f_{\leq N}$ for $P_{\leq N} f$ and similarly for the other operators.

The Littlewood-Paley operators commute with derivative operators, the free propagator, and complex conjugation.
They are self-adjoint and bounded on every $L^p$ and $\dot H^s_x$ space for $1 \leq p \leq \infty$ and $s\geq 0$.  They
also obey the following Sobolev and Bernstein estimates that we will use:
\begin{align*}
\|P_{\geq N} f\|_{L_x^p} &\lesssim N^{-s} \| |\nabla|^s P_{\geq N} f \|_{L_x^p},\\
\| |\nabla|^s  P_{\leq N} f\|_{L_x^p} &\lesssim N^{s} \| P_{\leq N} f\|_{L_x^p},\\
\| |\nabla|^{\pm s} P_N f\|_{L_x^p} &\sim N^{\pm s} \| P_N f \|_{L_x^p},\\
\|P_{\leq N} f\|_{L_x^q} &\lesssim N^{\frac{n}{p}-\frac{n}{q}} \|P_{\leq N} f\|_{L_x^p},\\
\|P_N f\|_{L_x^q} &\lesssim N^{\frac{n}{p}-\frac{n}{q}} \| P_N f\|_{L_x^p},
\end{align*}
whenever $s \geq 0$ and $1 \leq p \leq q \leq \infty$.

Let us temporarily fix a spacetime slab $\ir$.
We have already introduced two important norms on this slab, $\dot W(\ir)$ and $\dot N^1(\ir)$.  Now we introduce a few more.
We say that a pair of exponents $(q,r)$ is Schr\"odinger-\emph{admissible} if $\tfrac{2}{q} + \tfrac{n}{r} = \frac{n}{2}$
and $2 \leq q,r \leq \infty$. We define the $\dot{S}^0(\ir)$
\emph{Strichartz norm} by
\begin{equation}\label{S^0}
\|u\|_{\dot{S}^0(\ir)} = \|u\|_{\dot S^0} := \sup \Bigl(\sum_N \| P_N u \|^2_{q,r}\Bigr)^{1/2}
\end{equation}
where the supremum is taken over all admissible pairs $(q,r)$.  This Besov-type formulation of the Strichartz norm will be convenient for us
later when we need to prove nonlinear estimates.  We also define the $\dot S^1 (\ir)$ \emph{Strichartz norm}
to be
$$
\|u\|_{\dot S^1 (\ir)} := \| \nabla u \|_{\dot{S}^0(\ir)}.
$$
We define the associated spaces $\dot S^0(\ir)$, $\dot S^1(\ir)$ as the closure of the test functions under these norms.

We observe the inequality
\begin{equation}\label{square sum}
\|\bigl(\sum_N |f_N|^2 \bigr)^{1/2}\|_{q,r} \leq \bigl(\sum_N \|f_N\|^2_{q,r}\bigr)^{1/2}
\end{equation}
for all $2 \leq q,r \leq \infty$ and arbitrary functions $f_N$, which one proves by interpolating between the trivial
cases $(q,r) = (2,2), (2,\infty), (\infty,2), (\infty,\infty)$. In particular, \eqref{square sum} holds for all
admissible exponents $(q,r)$.  Combining this with the Littlewood-Paley inequality, we find
\begin{align*}
\| u \|_{q,r}& \lesssim \|\bigl(\sum_N |P_N u|^2\bigr)^{1/2}\|_{q,r}\\
            & \lesssim \bigl(\sum_N \|P_N u \|^2_{q,r}\bigr)^{1/2}\\
            & \lesssim \| u \|_{\dot S^0 (\ir)},
\end{align*}
which in particular implies
$$
\|\nabla u \| _{\lqlr (\ir)} \lesssim \|u\|_{\dot S^1 (\ir)}.
$$
{}From this and Sobolev embedding, the $\dot{S}^1$-norm controls the following spacetime norms:

\begin{lemma}\label{lemma strichartz norms}
For any $\dot{S}^1$ function $u$ on $\ir$, we have
\begin{align*}\label{strichartz norms}
\|\nabla u\|_{\infty,2} + \|u\|_{\dot W} + \|\nabla u\|_{\frac{2(n+2)}{n},\frac{2(n+2)}{n}}
        + \|\nabla u\|_{2, \frac{2n}{n-2}} & \\
 {} + \|u\|_{\infty, \frac{2n}{n-2}} + \|u\|_{\frac{2(n+2)}{n-2},\frac{2(n+2)}{n-2}} &\lesssim \|u\|_{\dot{S}^1}
\end{align*}
where all spacetime norms are on $\ir$.
\end{lemma}


By Lemma~\ref{lemma strichartz norms} and the definition of the Strichartz spaces $\dot S^0$ and $\dot S^1$,
we see that $\dot S^0 \subseteq C^0_t L^2_x$ and similarly, $\dot S^1 \subseteq C^0_t \dot H^1_x$,
by the usual limiting arguments. In particular $ \dot S^0 \cap \dot S^1 \subseteq C^0_t H^1_x$.

Next, let us recall the Strichartz estimates:

\begin{lemma}\label{lemma linear strichartz}
Let $I$ be a compact time interval, $k=0, 1$, and let $u : \ir \rightarrow \C$ be an $\dot{S}^k$ solution to
the forced Schr\"odinger equation
\begin{equation*}
i u_t + \Delta u = \sum_{m=1}^M F_m
\end{equation*}
for some functions $F_1 ,\dots,F_M$.  Then on $\ir$ we have
$$
\|u\|_{\dot{S}^k} \lesssim \|u(t_0)\|_{\dot{H}^k(\R^n)} + \sum_{m=1}^M \|\nabla^k F_m \|_{q'_m, r'_m}
$$
for any time $t_0 \in I$ and any admissible exponents $(q_1,r_1),\dots,(q_m,r_m)$. As
usual, $p'$ denotes the dual exponent to $p$, that is, $1/p + 1/p' = 1$.
\end{lemma}

Note that from this lemma and Lemma \ref{lemma strichartz norms}, we have
\begin{equation}\label{lls-special}
\| u \|_{\dot W} \lesssim \|u\|_{\dot S^1} \lesssim \| u(t_0) \|_{\dot H^1_x} + \| iu_t + \Delta u  \|_{\dot N^1}
\end{equation}
and, in particular, replacing $u$ with $u - e^{i(t-t_0)\Delta} u(t_0)$,
\begin{equation}\label{lls-special-2}
\| u - e^{i(t-t_0)\Delta} u(t_0) \|_{\dot W} \lesssim \| iu_t + \Delta u  \|_{\dot N^1}.
\end{equation}
In a similar spirit we have
\begin{equation}\label{lls-special-0}
\| u \|_{\frac{2(n+2)}{n}, \frac{2(n+2)}{n}} \lesssim \|u\|_{\dot S^0(\ir)} \lesssim \| u(t_0) \|_{L_x^2} +
\| iu_t + \Delta u  \|_{\frac{2(n+2)}{n+4}, \frac{2(n+2)}{n+4}}.
\end{equation}

\begin{proof}
To prove Lemma \ref{lemma linear strichartz}, let us first make the following reductions. We note that it suffices to
take $M=1$, since the claim for general $M$ follows from Duhamel's formula and the triangle inequality. We can also take $k$
to be 0, since the estimate for $k= 1$ follows by applying $\nabla$ to both sides of the equation and
noting that $\nabla$ commutes with $i\partial_t+\Delta$. As the Littlewood-Paley operators also commute with
$i\partial_t+\Delta$, we have
$$
(i\partial_t+\Delta)P_{N}u=P_N F_1
$$
for all dyadic $N$. Applying the standard Strichartz estimates (see \cite{tao:keel}), we get
\begin{align}\label{P N Strichartz}
\|P_N u\|_{q,r}\lesssim \|P_Nu(t_0)\|_{L_x^2}+\|P_N F_1\|_{q'_1, r'_1}
\end{align}
for all admissible exponents $(q,r)$ and $(q_1,r_1)$. Squaring \eqref{P N Strichartz}, summing in $N$, using the definition
of the $\dot{S}^0$-norm and the Littlewood-Paley inequality together with the dual of \eqref{square sum}, we get the
claim.
\end{proof}

%
%
%
%

\section{Local well-posedness}\label{lwp-sec}

The goal of this section is to establish the preliminary local well-posedness theorem, i.e., Theorem~\ref{lwp}.
The material here is standard but we include it for completeness.  More precise control
on the continuous dependence on the data will be given in later sections.  Throughout this section the hypotheses
are as in Theorem \ref{lwp}. All spacetime norms will be over $\ir$ unless otherwise specified.

\subsection{Existence}
We consider first the question of local existence of $\dot{S}^1$ solutions to
\eqref{equation 1} which we address via an iterative procedure. We define the following iterates
\begin{align*}
u^{(0)}(t)&:=0, \\
u^{(1)}(t)&:= e^{i(t-t_0)\Delta}u_0,
\end{align*}
and for $m\geq 1$,
\begin{align}\label{recurrence}
u^{(m+1)}(t):=e^{i(t-t_0)\Delta}u_0-i\int_{t_0}^t e^{i(t-s)\Delta}f(u^{(m)}(s))ds.
\end{align}

Let us first remark that the free evolution being small implies that all the iterates are small. Indeed, by
\eqref{lls-special-2} and the triangle inequality, followed by \eqref{chain-game} and \eqref{small}, we have
\begin{align*}
\|u^{(m+1)}\|_{\dot W}
&\lesssim \|e^{i(t-t_0)\Delta}u_0\|_{\dot W}+ \|f(u^{(m)})\|_{\dot N^1}\\
&\lesssim \eta + \|u^{(m)}\|_{\dot W}^{\frac{n+2}{n-2}}.
\end{align*}
Using \eqref{small} as the base case of an induction hypothesis and choosing $\eta_0$ sufficiently small, we deduce
\begin{align}\label{iterates small}
\|u^{(m)}\|_{\dot W}\lesssim \eta, \quad \|f(u^{(m)})\|_{\dot N^1} \lesssim \eta^{\frac{n+2}{n-2}}
\end{align}
for all $m\geq 1$.  Applying \eqref{lls-special}, we conclude that
\begin{equation}\label{its1}
\|u^{(m)}\|_{\dot{S}^1} \lesssim \| u_0 \|_{\dot H^1(\R^n)} + \eta^{\frac{n+2}{n-2}}
\end{equation}
for all $m \geq 1$.

Next, we consider differences of the form $u^{(m+1)}-u^{(m)}$. By the recurrence relation \eqref{recurrence}, in order to
estimate $u^{(m+1)}-u^{(m)}$ we need to control $f(u^{(m)})-f(u^{(m-1)})$.
By \eqref{iterates small}, \eqref{f diff} (with $u=u^{(m)}$ and $v=u^{(m-1)}$), and \eqref{lls-special-0}, we estimate
\begin{align*}
\|u^{(m+1)}-u^{(m)}\|_{\dot{S}^0}
&\lesssim \|f(u^{(m)})-f(u^{(m-1)})\|_{\frac{2(n+2)}{n+4},\frac{2(n+2)}{n+4}}\\
&\lesssim \|u^{(m)}-u^{(m-1)}\|_{\frac{2(n+2)}{n},\frac{2(n+2)}{n}}
\sum_{m' = m,m-1} \|u^{(m')}\|_{\frac{2(n+2)}{n-2},\frac{2(n+2)}{n-2}}^{\frac{4}{n-2}}\\
&\lesssim \|u^{(m)}-u^{(m-1)}\|_{\dot{S}^0}\sum_{m'=m,m-1} \|u^{(m')}\|_{\frac{2(n+2)}{n-2},\frac{2(n+2)}{n-2}}^{\frac{4}{n-2}} \\
&\lesssim \eta^{\frac{4}{n-2}}\|u^{(m)}-u^{(m-1)}\|_{\dot{S}^0},
\end{align*}
for all $m\geq 1$. Also, as $u_0 \in L^2(\R^n)$ by hypothesis, we have
$$ \| u^{(1)} \|_{\dot S^0} \lesssim \|u_0 \|_{L^2(\R^n)}$$
thanks to Lemma \ref{lemma linear strichartz}.
Choosing $\eta$ sufficiently small, we conclude that the sequence
$\{u^{(m)}\}_m$ is Cauchy in $\dot{S}^0(\ir)$ and hence, there exists $u\in \dot{S}^0(\ir)$ such that $u^{(m)}$
converges to $u$ in $\dot{S}^0(\ir)$ (and thus also in $L_{t,x}^{2(n+2)/n}(\ir)$) as $m\to\infty$ and furthermore,
we have \eqref{bounded-s0}. The above argument also shows that $f(u^{(m)})$
converges to $f(u)$ in $L_{t,x}^{2(n+2)/(n+4)}(\ir)$.  Applying \eqref{lls-special-0} again
we conclude that $u$ solves \eqref{iteration}.  As the $u^{(m)}$ converge strongly to $u$ in $\dot S^0(\ir)$ and
remain bounded in $\dot W(\ir)$ and $\dot S^1(\ir)$ thanks to \eqref{iterates small} and \eqref{its1}, we see that
the $u^{(m)}$ converge weakly to $u$ in $\dot W(\ir)$ and $\dot S^1(\ir)$ and we obtain the bounds
\eqref{small solution in L^p}, \eqref{bounded-s1dot}.  As $u$ lies in both $\dot S^1(\ir)$ and
$\dot S^0(\ir)$, it is in $C^0_t H^1_x(\ir)$ and is thus a strong solution to \eqref{equation 1}.

\subsection{Uniqueness} To prove uniqueness, we will in fact prove the following stronger (and standard) result:

\begin{proposition}[Unconditional uniqueness in $H^1(\R^n)$, \cite{cazenave:book}]\label{unconditional}
Let $I$ be a time interval containing $t_0$, and let $u_1, u_2 \in C^0_t H^1_x(\ir)$ be two strong solutions to
\eqref{equation 1} in the sense of \eqref{iteration}.  Then $u_1=u_2$ almost everywhere.
\end{proposition}

\begin{proof} By a standard continuity argument we may shrink $I$ so that one can construct the solution given by
the preceding iteration argument. Without loss of generality, we may to take $u_2$ to be this solution;
thus, $u_2 \in \dot S^1(\ir)$.

Write $v := u_1 - u_2 \in C^0_t H^1_x(\ir)$, and let $\Omega \subset I$ be the set of times $t$ for which
$v(t) = 0$ almost everywhere.  Then $\Omega$ is closed and contains $t_0$, so it suffices to show that
$\Omega$ is open in $I$.  Let $t_1 \in \Omega$ and suppose also that $[t_1,t_1+\tau) \subset I$ for some small $\tau > 0$.
Henceforth, we will work entirely on the slab $[t_1,t_1+\tau) \times \R^n$.
By shrinking $\tau$ as much as necessary, we may use the continuity and vanishing of $v$ at $t_1$ to assume
$$ \| v \|_{L^\infty_t \dot H^1_x} \leq \eta$$
where $\eta$ is a small absolute constant to be chosen later; in particular, from Sobolev embedding we have
$$ \| v \|_{\infty, \frac{2n}{n-2}} \leq \eta.$$
Similarly, from the hypothesis $u_2 \in \dot S^1(\ir)$ and
Lemma \ref{lemma strichartz norms} and \eqref{bounded-s1dot} we may take
$$ \| u_2 \|_{\frac{2(n+2)}{n-2}, \frac{2(n+2)}{n-2}} \lesssim \eta.$$

{}From \eqref{iteration} and the hypothesis $u_1, u_2 \in C^0_t H^1_x$, we easily see that
$$ u_j(t) = e^{i(t-t_1)} u_j(t_1) -i \int_{t_1}^t e^{i(t-s)\Delta} f(u_j(s))\ ds$$
for $j=1,2$.  Subtracting these and recalling that $v(t_1) = 0$, we conclude that
$$ v(t) = -i \int_{t_1}^t e^{i(t-s)\Delta} (f(u_1(s)) - f(u_2(s)))\ ds.$$
{}From \eqref{f diff} we have
$$ f(u_1(s)) - f(u_2(s)) = O( |v(s)|^{(n+2)/(n-2)} ) + O( |u_2(s)|^{4/(n-2)} |v(s)| ).$$
Applying Lemmas \ref{lemma strichartz norms} and \ref{lemma linear strichartz} followed by H\"older's inequality,
we conclude that
\begin{align*}
\|v\|_{2, \frac{2n}{n-2}}
&\lesssim \| |v|^{(n+2)/(n-2)} \|_{2, \frac{2n}{n+2}}
+ \| |u_2|^{4/(n-2)} |v| \|_{\frac{2(n+2)}{n+6}, \frac{2n(n+2)}{n^2+4n-4}} \\
&\lesssim \| v \|_{2, \frac{2n}{n-2}}
( \|v\|_{\infty, \frac{2n}{n-2}} + \|u_2\|_{\frac{2(n+2)}{n-2}, \frac{2(n+2)}{n-2}} )^{4/(n-2)}\\
&\lesssim \eta^{4/(n-2)} \| v \|_{2, \frac{2n}{n-2}}.
\end{align*}
As $v$ was already finite in $L^\infty_t L^{2n/(n-2)}_x$ (and hence in $L^2_t L^{2n/(n-2)}_x$), we conclude (by taking $\eta$ small enough)
that $v$ vanishes almost everywhere
on $[t_1,t_1+\tau)$.  This shows that $\Omega$ is open in the forward direction. A similar argument establishes the
openness in the backwards direction, thus concluding the proof.
\end{proof}

\begin{remark}
By combining this uniqueness statement with the existence theorem, one can show that for the Cauchy problem \eqref{equation 1} there is a unique
maximal interval $I$ containing $t_0$, such that the slab $\ir$ supports a strong solution.  Furthermore, $I$ is open and the solution
has finite $\dot W(J \times \R^n)$-norm for any  compact $J \subset I$.  It is also true that if $I$ has a finite endpoint, then
the $\dot W$-norm will blow up near that endpoint (see Lemma \ref{blow} below).
\end{remark}

\subsection{Continuous dependence in rough norms}

We turn now to the Lipschitz bound \eqref{lipschitz}. Again we write $v := u - \tilde u$.  By \eqref{f diff},
\eqref{lls-special-0}, and H\"older's inequality, we have
\begin{align*}
\|v\|_{\dot{S}^0}&\lesssim \|u_0-\tilde u_0\|_{L_x^2} + \|f(u)-f(\tilde u)\|_{\frac{2(n+2)}{n+4}, \frac{2(n+2)}{n+4}}\\
&\lesssim \|u_0-\tilde u_0\|_{L_x^2}+ \bigl(\|u\|_{\frac{2(n+2)}{n-2},\frac{2(n+2)}{n-2}}^{\frac{4}{n-2}}
        +\|\tilde u\|_{\frac{2(n+2)}{n-2},\frac{2(n+2)}{n-2}}^{\frac{4}{n-2}}\bigr) \|v\|_{\frac{2(n+2)}{n},\frac{2(n+2)}{n}}\\
&\lesssim \|u_0-\tilde u_0\|_{L_x^2} + \bigl(\|u\|_{\frac{2(n+2)}{n-2},\frac{2(n+2)}{n-2}}^{\frac{4}{n-2}}
        +\|\tilde u\|_{\frac{2(n+2)}{n-2},\frac{2(n+2)}{n-2}}^{\frac{4}{n-2}}\bigr) \|v\|_{\dot{S}^0}.
\end{align*}
Applying \eqref{sob-embed} and \eqref{small solution in L^p}, we conclude
$$ \|v\|_{\dot S^0} \lesssim \|u_0-\tilde u_0\|_{L_x^2} + \eta^{4/(n-2)} \|v\|_{\dot{S}^0}.$$
By taking $\eta_0$ small enough, we obtain \eqref{lipschitz} as desired.
This concludes the proof of Theorem~\ref{lwp}.
\hfill$\square$

We end this section with the following companion to Theorem \ref{lwp}.

\begin{lemma}[Standard blowup criterion, \cite{cw0}, \cite{cwI}]\label{blow}
Let $u_0\in H^1_x$ and let $u$ be a strong $\dot S^1$ solution to \eqref{equation 1} on the slab $[t_0, T_0]\times\R^n$ such that
\begin{align}\label{norm finite}
\|u\|_{L_{t,x}^{\frac{2(n+2)}{n-2}}([t_0, T_0]\times\R^n)} < \infty.
\end{align}
Then there exists $\delta=\delta(u_0)>0$ such that the solution $u$ extends to a strong $\dot S^1$ solution to
\eqref{equation 1} on the slab $[t_0, T_0+\delta]\times\R^n$.
\end{lemma}

In the contrapositive, this lemma asserts that if a solution cannot be continued strongly beyond a time $T_*$, then the
$L_{t,x}^{\frac{2(n+2)}{n-2}}$-norm must blow up at that time.  One can also establish that other scale-invariant norms
(except for those norms involving $L^\infty_t$) also blow up at this time, but we will not do so here.

\begin{proof}  Let us denote the norm in \eqref{norm finite} by $M$.
The first step is to establish an $\dot{S}^1$ bound on $u$. In order to do so, we subdivide $[t_0, T_0]$ into
$N\sim \bigl(1+\frac{M}{\nu}\bigr)^{\frac{2(n+2)}{n-2}}$ subintervals $J_k$ such that
\begin{align}\label{first}
\|u\|_{L_{t,x}^{\frac{2(n+2)}{n-2}}(J_k\times\R^n)}\leq \nu
\end{align}
where $\nu$ is a small positive constant.
By \eqref{chain-game} and \eqref{lls-special}, we have
\begin{align*}
\|u\|_{\dot{S}^1(J_k\times\R^n)}
&\lesssim \|u(t_k)\|_{\dot{H}^1(\R^n)} +\|f(u)\|_{\dot N^1(J_k\times\R^n)} \\
&\lesssim \|u(t_k)\|_{\dot{H}^1(\R^n)} +\|u\|_{L_{t,x}^{\frac{2(n+2)}{n-2}}(J_k\times\R^n)}^{\frac{4}{n-2}}\| u\|_{\dot S^1(J_k\times\R^n)}\\
&\lesssim \|u(t_k)\|_{\dot{H}^1(\R^n)} + \nu^{\frac{4}{n-2}}\|u\|_{\dot S^1(J_k\times\R^n)}
\end{align*}
for each interval $J_k$ and any $t_k \in J_k$.  If $\nu$ is sufficiently small, we conclude
$$
\|u\|_{\dot{S}^1(J_k\times\R^n)}\lesssim \|u(t_k)\|_{\dot{H}^1_x},
$$
Recall that the $\dot S^1$-norm controls the $L^\infty_t \dot H^1_x$-norm.  Thus, we may glue these bounds together inductively to obtain a bound of the form
$$
\|u\|_{\dot{S}^1([t_0, T_0]\times\R^n)}\leq C(\|u_0\|_{\dot{H}^1_x},M,\nu),
$$
which by Lemma~\ref{lemma strichartz norms} implies
\begin{align}\label{bound}
\|u\|_{\dot W([t_0, T_0]\times\R^n)}\leq C(\|u_0\|_{\dot{H}^1_x},M,\nu).
\end{align}

Now let $t_0\leq \tau<T_0$. By \eqref{chain-game} and \eqref{lls-special-2}, we have
\begin{align*}
\|u - e^{i(t-\tau)\Delta}u(\tau)\|_{\dot W([\tau, T_0]\times\R^n)}
&\lesssim \|f(u)\|_{\dot N^1([\tau, T_0]\times\R^n)}\\
&\lesssim \|u\|_{\dot W([\tau, T_0]\times\R^n)}^{\frac{n+2}{n-2}}
\end{align*}
and thus, by the triangle inequality,
$$ \| e^{i(t-\tau)\Delta}u(\tau)\|_{\dot W([\tau, T_0]\times\R^n)}
\lesssim \|u\|_{\dot W([\tau, T_0]\times\R^n)}^{\frac{n+2}{n-2}} + \|u\|_{\dot W([\tau, T_0]\times\R^n)}.$$
Let $\eta_0$ be as in Theorem~\ref{lwp}. By \eqref{bound}, taking $\tau$
sufficiently close to $T_0$, we obtain
$$
\|e^{i(t-\tau)\Delta}u_0\|_{\dot W([\tau, T_0]\times\R^n)}\leq \eta_0/2
$$
(say), while from Strichartz inequality we have
$$
\|e^{i(t-\tau)\Delta}u_0\|_{\dot W(\R \times \R^n)} < \infty.
$$
By the monotone convergence theorem, we deduce that there exists $\delta=\delta(u_0)>0$ such that
$$
\|e^{i(t-\tau)\Delta}u_0\|_{\dot W([\tau, T_0+\delta]\times\R^n)}\leq \eta_0.
$$
By Theorem~\ref{lwp}, there exists a unique solution to \eqref{equation 1} with initial data $v(\tau)$ at time $t=\tau$
which belongs to $\dot{S}^1([\tau, T_0+\delta]\times\R^n)$. By Proposition \ref{unconditional}, we see that $u=v$ on
$[\tau, T_0]\times\R^n$ and thus $v$ is an extension of $u$ to $[t_0, T_0+\delta]\times\R^n$.
\end{proof}

%
%
%
%

\section{Short-time perturbations}\label{short-sec}

The goal of this section is to prove Theorem~\ref{short-time theorem}. By the well-posedness theory that we have
developed in the previous section, it suffices to prove \eqref{close in L^p-0}-\eqref{one}
as \emph{a priori} estimates, that is, we assume that the solution $u$ already exists and belongs to $\dot{S}^1(\ir)$.

\begin{remark}\label{redundant}
By \eqref{lls-special} and Plancherel's theorem we have
\begin{align*}
\Bigl(\sum_N \|P_N e^{i(t-t_0)\Delta}\bigl(u(t_0)-\util(t_0)&\bigr)\|^2_{\dot W(\ir)}\Bigr)^{1/2}\\
&\lesssim \Bigl(\sum_N \|P_N \nabla(u(t_0)-\util(t_0)\bigr)\|^2_{\infty,2}\Bigr)^{1/2}\\
&\lesssim \|\nabla (u(t_0)-\util(t_0)\bigr)\|_{\infty,2} \\
&\lesssim E'
\end{align*}
on the slab $\ir$,
so the hypothesis \eqref{closer-0} is redundant if $E'=O(\eps)$.
\end{remark}

By time symmetry, we may assume that $t_0=\inf I$. We will first give a simple proof of Theorem~\ref{short-time theorem}
in dimensions $3\leq n\leq 6$ (following the arguments in \cite{gopher}, \cite{rv} covering the cases $n=3,4$ respectively).
Let $v := u - \tilde u$. Then $v$ satisfies the following initial value problem:
\begin{equation}\label{equation diff}
\begin{cases}
i v_t +\Delta v = f(\util+v)-f(\util)-e \\
v(t_0,x) = u(t_0, x)-\util(t_0,x).
\end{cases}
\end{equation}
For $T \in I$ define
\begin{align*}
S(T) := \| (i \partial_t + \Delta)v + e \|_{\dot N^1([t_0,T]\times \R^n)}.
\end{align*}
We will now work entirely on the slab $[t_0,T] \times \R^n$.
By \eqref{closer-0}, \eqref{error small-0}, and \eqref{lls-special}, we get
\begin{align}\label{v bound}
\|v\|_{\dot W}
&\lesssim \|e^{i(t-t_0)\Delta }v(t_0)\|_{\dot W} + \|(i \partial_t +\Delta)v+ e\|_{\dot N^1} +
\|e\|_{\dot N^1} \\
&\lesssim S(T) + \eps,\notag
\end{align}
where we used \eqref{square sum} to estimate
\begin{equation}\label{free ss}
\begin{split}
\|e^{i(t-t_0)\Delta }v(t_0)\|_{\dot W}
&\lesssim \Bigl(\sum_N \|P_N e^{i(t-t_0)\Delta}\bigl(u(t_0)-\util(t_0)\bigr)\|^2_{\dot W}\Bigr)^{1/2}\\
&\lesssim \eps.
\end{split}
\end{equation}
By \eqref{sob-embed} and \eqref{v bound} we have
\begin{align}\label{v bound-1}
\| v\|_{\frac{2(n+2)}{n-2},\frac{2(n+2)}{n-2}} \lesssim S(T) + \eps.
\end{align}
On the other hand, from \eqref{chain} we have
\begin{align*}
\nabla\bigl[(i \partial_t + \Delta)v +e\bigr]
&=\nabla \bigl[f(\util + v) -f(\util)\bigr]\\
&= f_{z}(\util + v)\nabla(\util + v) +f_{\bar{z}}(\util +v)\nabla\overline{(\util + v)}\\
&\quad -f_{z}(\util)\nabla\util -f_{\bar{z}}(\util)\nabla\bar{\util},
\end{align*}
so, by our hypotheses on $f$, specifically \eqref{f_z} and \eqref{f_z diff}, we get
\begin{align}
\bigl|\nabla\bigl[(i \partial_t + \Delta)v +e\bigr]\bigr|
&\lesssim |\nabla \util|\bigr( | f_{z}(\util + v)-f_{z}(\util)|+|f_{\bar{z}}(\util +v)-f_{\bar{z}}(\util)|\bigr)\notag\\
&\quad+|\nabla v| \bigl(|f_{z}(\util + v)|+|f_{\bar{z}}(\util +v)|\bigr)\notag\\
&\lesssim |\nabla \util||v|^{\frac{4}{n-2}} +|\nabla v||\util+v|^{\frac{4}{n-2}}\label{v!}.
\end{align}
Hence by \eqref{holder}, \eqref{finite S norm-0}, \eqref{v bound}, and \eqref{v bound-1}, we estimate
\begin{align*}
S(T)&\lesssim \|\util\|_{\dot W} \|v\|^{\frac{4}{n-2}}_{\frac{2(n+2)}{n-2},\frac{2(n+2)}{n-2}}
+\|v\|_{\dot W} \|\util\|^{\frac{4}{n-2}}_{\frac{2(n+2)}{n-2},\frac{2(n+2)}{n-2}}+
\|v\|_{\dot W} \|v\|^{\frac{4}{n-2}}_{\frac{2(n+2)}{n-2},\frac{2(n+2)}{n-2}}\\
&\lesssim \eps_0(S(T)+\eps)^{\frac{4}{n-2}} +\eps_0^{\frac{4}{n-2}}(S(T)+\eps)+(S(T)+\eps)^{\frac{n+2}{n-2}}.
\end{align*}
If $\frac{4}{n-2}\geq 1$, i.e., $3 \leq n\leq 6$, a standard continuity argument shows that if we take
$\eps_0 = \eps_0(E,E')$ sufficiently small we obtain
\begin{align}\label{S bound}
S(T) \leq \eps \quad \text{for all} \ T \in I,
\end{align}
which implies \eqref{one}.  Using \eqref{v bound} and \eqref{S bound},
one easily derives \eqref{close in L^p-0}. To obtain \eqref{close in S^1-0}, we use \eqref{close-0}, \eqref{error small-0},
\eqref{lls-special}, and \eqref{S bound}:
\begin{align*}
\|u-\tilde{u}\|_{\dot{S}^1(\ir)}
&\lesssim \|u(t_0)-\tilde{u}(t_0)\|_{\dot{H}^1_x}+\bigl\|(i \partial_t +\Delta)v+ e\bigr\|_{\dot N^1(\ir)}+\|e\|_{\dot N^1(\ir)}\\
&\lesssim E'+S(t)+\eps\\
&\lesssim E'+\eps.
\end{align*}
By the triangle inequality, \eqref{sob-embed}, \eqref{finite S norm-0}, and \eqref{v bound}, we have
$$
\| u\|_{L_{t,x}^{\frac{2(n+2)}{n-2}}(\ir)} \lesssim \|u\|_{\dot W(\ir)}\lesssim \eps+\eps_0.
$$
Another application of \eqref{chain-game} and \eqref{lls-special}, as well as \eqref{finite energy-0}, \eqref{close-0} yields
\begin{align*}
\|u \|_{\dot{S}^1(\ir)}
&\lesssim \| u (t_0)\|_{\dot H^1_x} + \|f(u)\|_{\dot N^1(\ir)}\\
&\lesssim E+E' + \|u\|_{\dot W(\ir)}^{(n+2)/(n-2)}\\
& \lesssim E+E'+(\eps+\eps_0)^{\frac{n+2}{n-2}},
\end{align*}
which proves \eqref{u in S^1-0}, provided $\eps_0$ is sufficiently small depending on $E$ and $E'$.

This concludes the proof of Theorem~\ref{short-time theorem} in dimensions $3\leq n\leq 6$. In order to prove the
theorem in higher dimensions, we are forced to avoid taking a full derivative since this is what turns the nonlinearity from
Lipschitz into just H\"older continuous of order $\frac{4}{n-2}$.  Instead, we must take fewer than $\frac{4}{n-2}$ derivatives.
As we still need to iterate in spaces that scale like $\dot{S}^1$, we either have to increase the space or the time integrability of the
usual Strichartz norms. The option of increasing the spatial integrability is suggested by the exotic Strichartz estimates of
Foschi, \cite{foschi}, but it turns out to be somewhat easier to increase the time integrability instead; this idea was used in the closely
related context of the energy-critical non-linear Klein-Gordon equation by Nakanishi \cite{nakanishi}.  We will choose the norm
$X = X(\ir)$ defined by
$$
\|u\|_{X} := \Bigl(\sum_N N^{8/(n+2)} \| P_N u \|_{n+2, \frac{2(n+2)}{n}}^2 \Bigr)^{1/2}.
$$
We observe that this norm is controlled by the $\dot S^1$-norm.
Indeed, by Sobolev embedding, the boundedness of the Riesz transforms on every
$L_x^p$, $1<p<\infty$, the dual of \eqref{square sum}, and the definition of the $\dot{S}^0$-norm, we get
\begin{align*}
\| u \|_{X}
& = \Bigl(\sum_N N^{8/(n+2)} \| P_N u \|_{n+2,\frac{2(n+2)}{n}}^2 \Bigr)^{1/2}\\
&\sim \Bigl(\sum_N \| |\nabla|^{4/(n+2)} P_N u \|_{n+2, \frac{2(n+2)}{n}}^2 \Bigr)^{1/2}\\
&\lesssim \Bigl(\sum_N \| \nabla P_N u \|_{n+2, \frac{2n(n+2)}{n^2+2n+4}}^2 \Bigr)^{1/2}\\
&\lesssim \|u \|_{\dot{S}^1(\ir)}
\end{align*}
since $(n+2, \frac{2n(n+2)}{n^2+2n+4})$ is a Schr\"odinger-admissible pair.

Now we need an inhomogeneous estimate.  We use the norm
$$
\|F\|_{Y} := \Bigl(\sum_N N^{8/(n+2)} \| P_N F \|_{\frac{n+2}{3}, \frac{2(n+2)}{n+4}}^2\Bigr)^{1/2}.
$$
Just as the $X$-norm is a variant of the $\dot S^1$-norm, one should think of the $Y$-norm as a variant of the
$\dot N^1$-norm.  The reason we use
these norms instead of the usual $\dot S^1$ and $\dot N^1$ norms is that they require roughly $4/(n+2) < 4/(n-2)$
degrees of differentiability only,
while still having the same scaling as the full-derivative spaces $\dot S^1$ and $\dot N^1$.  This
will become relevant when we have to address the limited regularity available for $f_z$ and $f_{\bar{z}}$.

\begin{lemma}[Exotic Strichartz estimate]  For any $F \in Y$, we have
\begin{align}\label{inhom Strichartz}
\Bigl\|\int_{t_0}^t e^{i(t-s)\Delta} F(s)ds\Bigr\|_{X}\lesssim \|F\|_{Y}
\end{align}
\end{lemma}

\begin{proof}
Interpolating between \eqref{dispersive ineq} and the conservation of mass,
we get the following dispersive inequality
$$
\| e^{i(t-s)\Delta} F(s) \|_{\frac{2(n+2)}{n}} \lesssim |t-s|^{-\frac{n}{n+2}} \|F(s)\|_{\frac{2(n+2)}{n+4}}
$$
whenever $t \neq s$.
By fractional integration, this implies that on the slab $\ir$ we have
$$
\Bigl\| \int_{t_0}^t e^{i(t-s)\Delta} F(s)\ ds \Bigr\|_{n+2, \frac{2(n+2)}{n}} \lesssim \| F \|_{\frac{n+2}{3}, \frac{2(n+2)}{n+4}}.
$$
As the Littlewood-Paley operators commute with the free evolution, we get
$$
\Bigl\| P_N \int_{t_0}^t  e^{i(t-s)\Delta} F(s)\ ds \Bigr\|_{n+2, \frac{2(n+2)}{n}} \lesssim \| P_N F \|_{\frac{n+2}{3}, \frac{2(n+2)}{n+4}}.
$$
Squaring the above inequality, multiplying by $N^{8/(n+2)}$, and summing over all dyadic $N$'s, we obtain
\eqref{inhom Strichartz}.
\end{proof}

To complete the set of tools needed to prove Theorem \ref{short-time theorem} in high dimensions, we need the following fractional variant of
\eqref{holder}.

\begin{lemma}[Nonlinear estimate]  On any slab $\ir$, we have
\begin{equation}\label{XY0}
\|f_z(v) u\|_{Y} \lesssim \| v \|_{\dot W(\ir)}^{\frac{4}{n-2}} \| u \|_{X}.
\end{equation}
whenever the right-hand side makes sense.  A similar statement holds with $f_z$ replaced by $f_{\bar{z}}$.
\end{lemma}

\begin{proof}  We just prove the claim for $f_z$, as the corresponding claim for $f_{\bar{z}}$ is identical.
The main estimate is

\begin{lemma}[Frequency-localized nonlinear estimate] For any dyadic $M$, we have
\begin{equation}\label{para}
\| P_M \bigl(f_z(v) u\bigr)\|_{\frac{n+2}{3}, \frac{2(n+2)}{n+4}}
\lesssim \|v\|_{\dot W}^{\frac{4}{n-2}}
 \sum_N \min\bigl(1, \bigl(\tfrac{N}{M}\bigr)^{\frac{4}{n-2}}\bigr) \| P_N u \|_{n+2, \frac{2(n+2)}{n}},
\end{equation}
where all spacetime norms are on $\ir$.
\end{lemma}

\begin{proof} We may rescale\footnote{Note that while the nonlinearity $f$ is not necessarily invariant
under the scaling $f(z) \mapsto \lambda^{(n+2)/(n-2)} f(z/\lambda)$, the
\emph{hypotheses} on the nonlinearity are invariant under this scaling and so, the
scaling argument remains justified in such cases.} such that $M=1$. Splitting $u$ into
Littlewood-Paley components, we estimate

\begin{align}\label{split}
\| P_1 \bigl(f_z(v) u\bigr)\|_{\frac{n+2}{3}, \frac{2(n+2)}{n+4}}\lesssim \sum_N \| P_1 \bigl( f_z(v) P_N u \bigr) \|_{\frac{n+2}{3}, \frac{2(n+2)}{n+4}}.
\end{align}
Let us first consider the components where $N \geq 1/4$; we use \eqref{f_z}, H\"older's inequality, and \eqref{sob-embed} to bound the contribution
of these components to the right-hand side of \eqref{split} by
\begin{align*}
\sum_{N \geq 1/4} \| f_z(v) \|_{\frac{n+2}{2}, \frac{n+2}{2}} \| P_N u \|_{n+2, \frac{2(n+2)}{n}}
&\lesssim \sum_{N \geq 1/4} \| v \|_{\frac{2(n+2)}{n-2}, \frac{2(n+2)}{n-2}}^{4/(n-2)} \| P_N u \|_{n+2, \frac{2(n+2)}{n}} \\
&\lesssim \| v \|_{\dot W}^{\frac{4}{n-2}} \sum_{N\geq 1/4}\| P_N u \|_{n+2, \frac{2(n+2)}{n}}.
\end{align*}

Next, we consider the components where $N \leq 1/8$.  We can freely replace $f_z(v)$ by
$P_{\geq 1/2} f_z(v)$ and then drop the projection $P_1$, leaving us with the task of
estimating
$$
\sum_{N \leq 1/8} \| ( P_{\geq 1/2} f_z(v) ) P_N u \|_{\frac{n+2}{3}, \frac{2(n+2)}{n+4}}.
$$
We use H\"older to estimate this contribution by
$$ \sum_{N \leq 1/8} \| P_{\geq 1/2} f_z(v) \|_{\frac{n+2}{2}, p}
\| P_N u \|_{n+2, r}
$$
where $\frac{1}{r} :=\frac{n}{2(n+2)} - \frac{4}{n(n-2)}$ and $p := \frac{n(n-2)(n+2)}{2(n^2+4)}$.
{} From Bernstein, we have
$$
\| P_N u \|_{n+2, r} \lesssim  N^{\frac{4}{n-2}} \| P_N u \|_{n+2, \frac{2(n+2)}{n}},
$$
so the claim follows if we can establish that
\begin{align}\label{showshow}
\| P_{\geq 1/2} f_z(v) \|_{\frac{n+2}{2}, p}
\lesssim \|v\|_{\dot W}^{\frac{4}{n-2}}.
\end{align}
In order to prove \eqref{showshow}, we split $v$ into $v_{lo} := P_{<1} v$ and $v_{hi} := P_{\geq 1} v$, and use \eqref{f_z diff} to write
$$
f_z(v) = f_z(v_{lo})  + O( |v_{hi}|^{4/(n-2)} ).
$$
To deal with the contribution of the error term $O( |v_{hi}|^{4/(n-2)} )$, we discard $P_{\geq 1/2}$ and use H\"older and
Bernstein to estimate
\begin{align*}
\| P_{\geq 1/2} O( |v_{hi}|^{4/(n-2)} )\|_{\frac{n+2}{2}, p}
&\lesssim \| v_{hi} \|_{\frac{2(n+2)}{n-2}, \frac{2n(n+2)}{n^2+4}}^{\frac{4}{n-2}} \\
&\lesssim \| \nabla v_{hi} \|_{\frac{2(n+2)}{n-2}, \frac{2n(n+2)}{n^2+4}}^{\frac{4}{n-2}} \\
&\lesssim \| v \|_{\dot W}^{4/(n-2)}.
\end{align*}

We are thus left with the main term. We need to show
\begin{align}\label{show!}
\| P_{\geq 1/2} f_z(v_{lo}) \|_{\frac{n+2}{2}, p}
\lesssim \| f \|_{\dot W}^{\frac{4}{n-2}}.
\end{align}
In order to prove \eqref{show!}, we make a few remarks. First, observe that for any spatial function $F$ we have the bound
$$
\| P_{\geq 1/2} F \|_{p} \lesssim \sup_{|h| \leq 1} \| \tau_h F - F \|_{p}
$$
where $\tau_h F(x) := F(x-h)$ is the translation operator.  Indeed, from the triangle inequality we have
$$
\| \tau_y F - F \|_{p} \leq \langle y \rangle \sup_{|h| \leq 1} \| \tau_h F - F \|_{p}
$$
for any $y \in \R^n$; integrating this against the convolution kernel of $P_{<1/2}$ (which is rapidly decreasing and
has total mass one) we obtain the claim.  Next, observe from \eqref{f_z diff} that we have the pointwise
estimate
$$
\tau_h f_z(v_{lo}) - f_z(v_{lo}) = f_z(\tau_h v_{lo}) - f_z(v_{lo}) = O( |\tau_h v_{lo} -v_{lo}|^{\frac{4}{n-2}} )
$$
and thus, for any $t \in I$,
$$
\bigl\| P_{\geq 1/2} f_z(v_{lo}(t)) \bigr\|_{p}
\lesssim \sup_{|h| \leq 1} \| \tau_h v_{lo}(t)- v_{lo}(t) \|_{\frac{2n(n+2)}{n^2+4}}^{\frac{4}{n-2}}.
$$
But from the fundamental theorem of calculus and Minkowski's inequality, for $|h| \leq 1$ we have
$$
\| \tau_h v_{lo}(t)- v_{lo}(t) \|_{\frac{2n(n+2)}{n^2+4}} \leq \| \nabla v_{lo}(t) \|_{\frac{2n(n+2)}{n^2+4}}.
$$
Thus,
$$
\bigl\| P_{\geq 1/2} f_z(v_{lo}(t)) \bigr\|_{p}
\lesssim \| \nabla v_{lo}(t) \|_{\frac{2n(n+2)}{n^2+4}}^{\frac{4}{n-2}}
$$
and hence, by H\"older's inequality in the time variable, \eqref{show!} follows:
\begin{align*}
\bigl\| P_{\geq 1/2} f_z(v_{lo}) \bigr\|_{\frac{n+2}{2}, p}
&\lesssim \| \nabla v_{lo} \|_{\frac{2(n+2)}{n-2}, \frac{2n(n+2)}{n^2+4}}^{\frac{4}{n-2}}\\
&\lesssim \| v \|_{\dot W}^{\frac{4}{n-2}}.
\end{align*}
\end{proof}

We now return to the proof of \eqref{XY0}, and rewrite \eqref{para} as
\begin{align*}
M^{4/(n+2)} \| &P_M \bigl( f_z(v) u ) \|_{\frac{n+2}{3}, \frac{2(n+2)}{n+4}}\\
&\lesssim \| v \|_{\dot W}^{\frac{4}{n-2}}
 \sum_N \min\bigl((\tfrac{M}{N})^{\frac{4}{n+2}},(\tfrac{N}{M})^{\frac{16}{n^2-4}}\bigr)  N^{4/(n+2)} \| P_N u \|_{n+2, \frac{2(n+2)}{n}}.
\end{align*}
The claim \eqref{XY0} then follows from Schur's test.
\end{proof}

We are now ready to resume the proof of Theorem~\ref{short-time theorem} in dimensions $n>6$. Recall that $v:=u-\util$
satisfies the initial value problem \eqref{equation diff} and hence,
$$
v(t)=e^{i(t-t_0)\Delta}v(t_0)-i\int_{t_0}^t e^{i(t-s)\Delta}\bigl(f(\util+v)-f(\util)\bigr)(s)ds-i\int_{t_0}^t e^{i(t-s)\Delta}e(s)ds.
$$
We estimate
\begin{align*}
\|v\|_X
&\lesssim \|e^{i(t-t_0)\Delta}v(t_0)\|_X + \Bigl\|\int_{t_0}^t e^{i(t-s)\Delta}\bigl(f(\util+v)-f(\util)\bigr)(s)ds\Bigr\|_X\\
&\quad + \Bigl\|\int_{t_0}^t e^{i(t-s)\Delta}e(s)ds\Bigr\|_X,
\end{align*}
which by \eqref{inhom Strichartz} becomes
\begin{align}\label{generalized Strichartz}
\|v\|_X
\lesssim \|e^{i(t-t_0)\Delta}v(t_0)\|_X + \|f(\util+v)-f(\util)\|_Y +\Bigl\|\int_{t_0}^t e^{i(t-s)\Delta}e(s)ds\Bigr\|_X.
\end{align}

We consider first the free evolution term in \eqref{generalized Strichartz}. Using Sobolev embedding and the boundedness
of the Riesz transforms on $L_x^p$ for $1<p<\infty$, we estimate
\begin{align*}
\|e^{i(t-t_0)\Delta}v(t_0)\|_X
&=\Bigl(\sum_N N^{\frac{8}{n+2}} \| P_N e^{i(t-t_0)\Delta}v(t_0) \|_{n+2, \frac{2(n+2)}{n}}^2 \Bigr)^{1/2}\\
&\sim \Bigl(\sum_N  \| |\nabla| ^{\frac{4}{n+2}}P_N e^{i(t-t_0)\Delta}v(t_0) \|_{n+2, \frac{2(n+2)}{n}}^2 \Bigr)^{1/2}\\
&\lesssim \Bigl(\sum_N  \|  P_N \nabla e^{i(t-t_0)\Delta}v(t_0) \|_{n+2, \frac{2n(n+2)}{n^2+2n-4}}^2 \Bigr)^{1/2}.
\end{align*}
Now, we observe that $L^{n+2}_t L^{\frac{2n(n+2)}{n^2+2n-4}}_x$ interpolates between
$L^{\frac{2(n+2)}{n-2}}_t L^{\frac{2n(n+2)}{n^2+4}}_x$ and $L_t^\infty L_x^2$ and hence,
\begin{align*}
\|P_N \nabla & e^{i(t-t_0)\Delta}v(t_0)\|_{n+2, \frac{2n(n+2)}{n^2+2n-4}}\\
&\lesssim \|P_N \nabla e^{i(t-t_0)\Delta}v(t_0) \|_{\frac{2(n+2)}{n-2}, \frac{2n(n+2)}{n^2+4}}^{\frac{2}{n-2}}
   \|P_N \nabla e^{i(t-t_0)\Delta}v(t_0) \|_{\infty, 2}^{\frac{n-4}{n-2}}.
\end{align*}
The first factor is just the $\dot W(\ir)$ norm of $P_N e^{i(t-t_0)\Delta}v(t_0)$.
Squaring the above inequality, summing over all dyadic $N$'s, and applying H\"older's inequality for sequences, we obtain
\begin{align*}
\Bigl(\sum_N  \|  P_N \nabla  e^{i(t-t_0)\Delta}v(t_0)& \|_{n+2, \frac{2n(n+2)}{n^2+2n-4}}^2 \Bigr)^{1/2}\\
&\lesssim \Bigl(\sum_N  \|  P_N e^{i(t-t_0)\Delta}v(t_0) \|_{\dot W(\ir)}^2 \Bigr)^{\frac{2}{2(n-2)}}\\
&\qquad \times \Bigl(\sum_N  \|  P_N \nabla e^{i(t-t_0)\Delta}v(t_0) \|_{\infty, 2}^2 \Bigr)^{\frac{n-4}{2(n-2)}}.
\end{align*}
By \eqref{closer-0}, we have
$$
\Bigl(\sum_N  \|  P_N e^{i(t-t_0)\Delta}v(t_0) \|_{\dot W(\ir)}^2 \Bigr)^{\frac{2}{2(n-2)}}
\lesssim \eps^{\frac{2}{n-2}},
$$
while by the usual Strichartz estimates and \eqref{close-0}, we get
\begin{align*}
\Bigl(\sum_N  \|  P_N \nabla e^{i(t-t_0)\Delta}v(t_0) \|_{\infty, 2}^2 \Bigr)^{\frac{n-4}{2(n-2)}}
&\lesssim \Bigl(\sum_N  \|  P_N \nabla v(t_0) \|_{2}^2 \Bigr)^{\frac{n-4}{2(n-2)}}\\
&\lesssim \|\nabla v(t_0)\|_{2}^{\frac{n-4}{n-2}}\lesssim (E')^{\frac{n-4}{n-2}}.
\end{align*}
Hence,
\begin{align}\label{free guy}
\Bigl(\sum_N  \|  P_N \nabla  e^{i(t-t_0)\Delta}v(t_0) \|_{n+2, \frac{2n(n+2)}{n^2+2n-4}}^2 \Bigr)^{1/2}
\lesssim \eps^{\frac{2}{n-2}} (E')^{\frac{n-4}{n-2}}.
\end{align}

We consider next the error term in \eqref{generalized Strichartz} which we estimate via Sobolev embedding and the usual
Strichartz estimates, recalling that $\nabla$ commutes with the free propagator and that the Riesz transforms are bounded
on $L^p_x$ for every $1<p<\infty$.
\begin{align*}
\Bigl\|\int_{t_0}^t e^{i(t-s)\Delta}e(s)ds\Bigr\|_X
&\sim \Bigl(\sum_N  \Bigl\| |\nabla|^{\frac{4}{n+2}} \int_{t_0}^t e^{i(t-s)\Delta}P_N e(s)ds \Bigr\|_{n+2, \frac{2(n+2)}{n}}^2 \Bigr)^{1/2}\\
&\lesssim \Bigl(\sum_N  \Bigl\| \nabla \int_{t_0}^t e^{i(t-s)\Delta} P_N e(s)ds \Bigr\|_{n+2, \frac{2n(n+2)}{n^2+2n-4}}^2 \Bigr)^{1/2}\\
&\lesssim \Bigl(\sum_N  \| P_N e\|_{\dot N^1(\ir)}^2 \Bigr)^{1/2}.
\end{align*}
By \eqref{error small-0} and the dual of \eqref{square sum},
$$
\Bigl(\sum_N  \| P_N e\|_{\dot N^1(\ir)}^2 \Bigr)^{1/2}
\lesssim \|e\|_{\dot N^1(\ir)}
\lesssim \eps,
$$
so
\begin{align}\label{error term guy}
\Bigl\|\int_{t_0}^t e^{i(t-s)\Delta}e(s)ds\Bigr\|_X \lesssim \eps.
\end{align}

We turn now to the remaining term on the right-hand side of \eqref{generalized Strichartz}. From \eqref{ftoc} we have
$$
f(\util+v)-f(\util)=\int_0^1 \bigl[f_z(\util+\theta v)v+f_{\bar{z}}(\util+\theta v)\bar{v}\bigr] d\theta
$$
and so, by using Minkowski's inequality, \eqref{f_z}, and \eqref{XY0}, we get
\begin{align*}
\|f(\util+v)-f(\util)\|_Y
&\lesssim \Bigl(\|\util\|_{\dot W(\ir)}^{\frac{4}{n-2}}+\|u\|_{\dot W(\ir)}^{\frac{4}{n-2}}\Bigr)
  \|v\|_X.
\end{align*}
In order to bound $\|u\|_{\dot W(\ir)}$, we will first estimate
the free evolution of $u(t_0)$, i.e., $e^{i(t-t_0)\Delta}u(t_0)$. By the triangle inequality and \eqref{free ss}, we have
\begin{align*}
\|e^{i(t-t_0)\Delta}u(t_0)\|_{\dot W(\ir)}
&\lesssim \|e^{i(t-t_0)\Delta}v(t_0)\|_{\dot W(\ir)}
 +\|e^{i(t-t_0)\Delta}\util(t_0)\|_{\dot W(\ir)}\\
&\lesssim \eps+\|e^{i(t-t_0)\Delta}\util(t_0)\|_{\dot W(\ir)}.
\end{align*}
On the other hand, by \eqref{chain-game}, \eqref{finite S norm-0}, \eqref{error small-0}, and \eqref{lls-special-2},
we get
\begin{align*}
\|e^{i(t-t_0)\Delta}\util(t_0)\|_{\dot W(\ir)}
&\lesssim \|\util\|_{\dot W(\ir)}
 +\|f(\util)\|_{\dot N^1(\ir)}+ \| e\|_{\dot N^1(\ir)}\\
&\lesssim  \eps_0+\|\util\|_{\dot W(\ir)}^{\frac{n+2}{n-2}} + \eps\\
&\lesssim \eps_0,
\end{align*}
for $\eps_0$ sufficiently small. Hence
$$
\|e^{i(t-t_0)\Delta}u(t_0)\|_{\dot W(\ir)}
\lesssim \eps+ \eps_0\lesssim \eps_0.
$$

Assuming $\eps_0$ is sufficiently small depending on $E$ and $E'$, the hypotheses of Theorem~\ref{lwp} hold and hence its
conclusions hold as well; in particular, \eqref{small solution in L^p} holds, i.e.,
\begin{align}\label{u small again}
\|u\|_{\dot W(\ir)}\lesssim \eps_0.
\end{align}
Returning to our previous computations, by \eqref{finite S norm-0} and \eqref{u small again}, we get
\begin{align}\label{nonlinear guy}
\|f(\util+v)-f(\util)\|_Y\lesssim \eps_0^{\frac{4}{n-2}}\|v\|_X.
\end{align}

Considering \eqref{free guy}, \eqref{error term guy}, and \eqref{nonlinear guy}, \eqref{generalized Strichartz} becomes
\begin{align*}
\|v\|_X\lesssim \eps^{\frac{2}{n-2}}(E')^{\frac{n-4}{n-2}}+\eps+\eps_0^{\frac{4}{n-2}}\|v\|_X.
\end{align*}
Assuming that $\eps_0$ is sufficiently small depending on $E'$, a standard continuity argument yields
\begin{align}\label{v in X}
\|v\|_X\lesssim \eps^{\frac{2}{n-2}}(E')^{\frac{n-4}{n-2}}.
\end{align}
By Sobolev embedding and \eqref{square sum}, we conclude
\begin{equation}\label{v in smthg}
\begin{split}
\| v \|_{n+2, \frac{2n(n+2)}{n^2-8}}
&\lesssim
\bigl\| |\nabla|^{\frac{4}{n+2}} v \bigr\|_{n+2, \frac{2(n+2)}{n}} \\
&\lesssim \Bigl(\sum_N  \bigl\| |\nabla| ^{\frac{4}{n+2}}P_N v \bigr\|_{n+2, \frac{2(n+2)}{n}}^2 \Bigr)^{1/2} \\
&\lesssim \|v\|_X \\
&\lesssim \eps^{\frac{2}{n-2}}(E')^{\frac{n-4}{n-2}}.
\end{split}
\end{equation}

We are now ready to upgrade our bounds on $v$; first, we will show \eqref{close in L^p-0}. Indeed, by Strichartz's
inequality, \eqref{error small-0}, and \eqref{free ss}, we have
\begin{align}
\|v\|_{\dot W(\ir)}
&\lesssim \|e^{i(t-t_0)\Delta}v(t_0)\|_{\dot W(\ir)}\notag\\
&\quad + \|f(\util+v)-f(\util)\|_{\dot N^1(\ir)} +\|e\|_{\dot N^1(\ir)}\notag\\
&\lesssim \eps+\|f(\util+v)-f(\util)\|_{\dot N^1(\ir)}\label{nabla v}.
\end{align}
By \eqref{v!}, we have
$$
\bigl|\nabla [f(\util+v)-f(\util)]\bigr|\lesssim |\nabla \util| |v|^{\frac{4}{n-2}} +|\nabla v| |u|^{\frac{4}{n-2}}
$$
and hence
\begin{align*}
\|f(\util+v)-f(\util)\|_{\dot N^1(\ir)}
&\lesssim \|\nabla \util\|_{a,b}\| v \|_{n+2, \frac{2n(n+2)}{n^2-8}}^{\frac{4}{n-2}}\\
&\quad+\|v\|_{\dot W(\ir)}\|u\|_{\frac{2(n+2)}{n-2}, \frac{2(n+2)}{n-2}}^{\frac{4}{n-2}}
\end{align*}
where $\frac{1}{a}:=\frac{1}{2}-\frac{4}{n^2-4}$ and $\frac{1}{b}:=\frac{n+2}{2n}-\frac{2(n^2-8)}{n(n^2-4)}$. Note that
the pair $(a,b)$ is a Schr\"odinger admissible pair and hence
$$
\|\nabla \util\|_{a,b}\lesssim \|\util\|_{\dot{S}^1(\ir)}.
$$
However, by \eqref{chain-game}, \eqref{lls-special}, and our hypotheses \eqref{finite energy-0}, \eqref{finite S norm-0},
and \eqref{error small-0}, we have
\begin{align*}
\|\util\|_{\dot{S}^1(\ir)}
&\lesssim \|\util(t_0)\|_{\dot{H}^1_x} + \|f(\util)\|_{\dot N^1(\ir)} +\|e\|_{\dot N^1(\ir)}\\
&\lesssim E+\|\util\|_{\dot W(\ir)}^{\frac{n+2}{n-2}} + \eps\\
&\lesssim E+\eps_0^{\frac{n+2}{n-2}} + \eps_0
\end{align*}
and hence, for $\eps_0$ sufficiently small depending on $E$,
\begin{align}\label{util S^1}
\|\util\|_{\dot{S}^1(\ir)}\lesssim E.
\end{align}

Returning to our previous computations and using \eqref{u small again} (combined with \eqref{sob-embed})
and \eqref{util S^1}, we have
\begin{align*}
\|f(\util+v)-f(\util)\|_{\dot N^1(\ir)}
&\lesssim \|\nabla \util\|_{a,b}\| v \|_{n+2, \frac{2n(n+2)}{n^2-8}}^{\frac{4}{n-2}}\\
&\quad+\|v\|_{\dot W(\ir)}\|u\|_{\frac{2(n+2)}{n-2}, \frac{2(n+2)}{n-2}}^{\frac{4}{n-2}}\\
&\lesssim E \| v \|_{n+2, \frac{2n(n+2)}{n^2-8}}^{\frac{4}{n-2}} + \eps_0^{\frac{4}{n-2}}\|v\|_{\dot W(\ir)},
\end{align*}
which by \eqref{v in smthg} yields
\begin{align}
\|f(\util+v)-f(\util)\|_{\dot N^1(\ir)}
&\lesssim \eps^{\frac{8}{(n-2)^2}}(E')^{\frac{4(n-4)}{(n-2)^2}}E + \eps_0^{\frac{4}{n-2}}\| v\|_{\dot W(\ir)}.\label{aha!}
\end{align}

Returning to \eqref{nabla v}, we find
\begin{align*}
\|v\|_{\dot W(\ir)}
&\lesssim \eps+\eps^{\frac{8}{(n-2)^2}}(E')^{\frac{4(n-4)}{(n-2)^2}}E + \eps_0^{\frac{4}{n-2}}\| v\|_{\dot W(\ir)}.
\end{align*}
Assuming $\eps_0$ is sufficiently small depending on $E$ and $E'$, a standard continuity argument yields
\eqref{close in L^p-0}, i.e.,
\begin{align*}
\|v\|_{\dot W(\ir)}\lesssim \eps+ \eps^{\frac{7}{(n-2)^2}}.
\end{align*}
By \eqref{close in L^p-0} and \eqref{aha!}, we get \eqref{one}. Indeed,
\begin{align*}
\bigl\|(i\partial_t+\Delta)(u-\util)+e\bigr\|_{\dot N^1(\ir)}
&\lesssim \|f(\util+v)-f(\util)\|_{\dot N^1(\ir)}\\
&\lesssim \eps^{\frac{8}{(n-2)^2}}(E')^{\frac{4(n-4)}{(n-2)^2}}E + \eps_0^{\frac{4}{n-2}}\bigl(\eps+ \eps^{\frac{7}{(n-2)^2}}\bigr)\\
&\lesssim \eps+ \eps^{\frac{7}{(n-2)^2}},
\end{align*}
provided $\eps_0=\eps_0(E,E')$ is sufficiently small.

To prove \eqref{close in S^1-0}, we use Strichartz's inequality, \eqref{close-0}, \eqref{error small-0} and
\eqref{close in L^p-0}:
\begin{align*}
\|v\|_{\dot{S}^1(\ir)}
&\lesssim \|v(t_0)\|_{\dot{H}^1_x}+\bigl\|(i\partial_t+\Delta)(u-\util)+e\bigr\|_{\dot N^1(\ir)}\\
&\quad + \|e\|_{\dot N^1(\ir)}\\
&\lesssim E' + \eps+ \eps^{\frac{7}{(n-2)^2}}.
\end{align*}

By triangle inequality, \eqref{close in S^1-0} and \eqref{util S^1} imply \eqref{u in S^1-0}, provided $\eps_0$ is
chosen sufficiently small depending on $E$ and $E'$.

%
%
%
%

\section{Long-time perturbations}\label{large-sec}

 The goal of this section is to prove Theorem \ref{long-time theorem} and to derive
scattering results as corollaries. We will prove Theorem~\ref{long-time theorem} under the additional assumption
\begin{align}
\|u(t_0)\|_{L_x^2}<\infty. \label{finite mass}
\end{align}
This additional assumption can be removed by the usual limiting argument: approximating $u(t_0)$ in $\dot{H}^1_x$ by
$\{u_n(t_0)\}_n\subset H^1_x$ and applying Theorem~\ref{long-time theorem} (under finite mass assumptions) with
$\util=u_m$, $e=0$, and $u=u_n$, we obtain that the sequence of solutions $\{u_n\}_n$ to \eqref{equation 1} with
initial data $\{u_n(t_0)\}_n$ is Cauchy in $\dot{S}^1(\ir)$ and thus convergent to an $\dot{S}^1$ solution $u$
to \eqref{equation 1} with initial data $u(t_0)$ at time $t=t_0$.

We will derive Theorem~\ref{long-time theorem} from Theorem~\ref{short-time theorem} by an iterative procedure.
First, we will assume without loss of generality that $t_0=\inf I$. Let $\eps_0=\eps_0(E, 2E')$ be as in
Theorem~\ref{short-time theorem}. Note that we need to replace $E'$ by the slightly larger $2E'$ as the $\dot{H}^1_x$-norm
of $u(t)-\util(t)$ may possibly grow in time.

The first step is to establish an $\dot{S}^1$ bound on $\util$. In order to do so, we subdivide $I$ into
$N_0\sim \bigl(1+\frac{M}{\eps_0}\bigr)^{\frac{2(n+2)}{n-2}}$ subintervals $J_k$ such that
\begin{align}\label{second}
\|\util\|_{L_{t,x}^{\frac{2(n+2)}{n-2}}(J_k\times\R^n)}\leq \eps_0.
\end{align}
By  \eqref{chain-game}, \eqref{lls-special}, \eqref{finite energy}, \eqref{error small}, and \eqref{second}, we estimate
\begin{align*}
\|\util\|_{\dot{S}^1(J_k\times\R^n)}
&\lesssim \|\util(t_0)\|_{\dot{H}^1_x} +\|f(\util)\|_{\dot N^1(J_k\times\R^n)} +\|e\|_{\dot N^1(J_k\times\R^n)}\\
&\lesssim E +\|\util\|_{L_{t,x}^{\frac{2(n+2)}{n-2}}(J_k\times\R^n)}^{\frac{4}{n-2}}\|\util\|_{\dot W(J_k\times\R^n)} +\eps\\
&\lesssim E + \eps_0^{\frac{4}{n-2}}\|\util\|_{\dot{S}^1(J_k\times\R^n)} +\eps.
\end{align*}
A standard continuity argument yields
$$
\|\util\|_{\dot{S}^1(J_k\times\R^n)}\lesssim E,
$$
provided $\eps_0$ is sufficiently small depending on $E$. Summing these bounds over all the intervals $J_k$, we obtain
$$
\|\util\|_{\dot{S}^1(\ir)}\leq C(E,M,\eps_0),
$$
which by Lemma~\ref{lemma strichartz norms} implies
$$
\|\util\|_{\dot W(\ir)}\leq C(E,M,\eps_0).
$$
We now subdivide $I$ into $N_1=C(E, M, \eps_0)$ subintervals $I_j=[t_j, t_{j+1}]$ such that
\begin{align}\label{third}
\|\util\|_{\dot W(I_j\times\R^n)}\leq\eps_0.
\end{align}
Choosing $\eps_1$ sufficiently small depending on $N_1$, $E$, and $E'$, we apply  Theorem~\ref{short-time theorem}
to obtain for each $j$ and all $0<\eps<\eps_1$,
\begin{align*}
\|u-\util\|_{\dot W(I_j\times\R^n)} &\leq C(j)\bigl(\eps+\eps^{\frac{7}{(n-2)^2}}\bigr)\\
\|u-\util\|_{\dot{S}^1(I_j\times\R^n)} &\leq C(j)\bigl(E'+\eps+ \eps^{\frac{7}{(n-2)^2}}\bigr)\\
\|u\|_{\dot{S}^1(I_j\times\R^n)} &\leq C(j)(E+E')\\
\bigl\| (i\partial_t+\Delta)(u-\util)+e \bigr\|_{\dot N^1(I_j\times\R^n)} &\leq C(j)\bigl(\eps+\eps^{\frac{7}{(n-2)^2}}\bigr),
\end{align*}
provided we can show that \eqref{close} and \eqref{closer} hold with $t_0$ replaced by $t_j$. We verify this using
an inductive argument. By \eqref{close}, \eqref{lls-special}, and the inductive hypothesis,
\begin{align*}
\|u(t_{j+1})-\util(t_{j+1})\|_{\dot{H}^1_x}
&\lesssim \|u(t_0)-\util(t_0)\|_{\dot{H}^1_x} + \|e\|_{\dot N^1([t_0, t_{j+1}]\times\R^n)}\\
&\quad +\bigl\|(i\partial_t+\Delta)(u-\util)+e\bigr\|_{\dot N^1([t_0, t_{j+1}]\times\R^n)}\\
&\lesssim E' +\eps +\sum_{k=0}^jC(k)\bigl(\eps+\eps^{\frac{7}{(n-2)^2}}\bigr).
\end{align*}
As the Littlewood-Paley operators commute with derivative operators and the free propagator, by Strichartz, we estimate
\begin{align*}
\|& P_N e^{i(t-t_{j+1})\Delta}\bigl(u(t_{j+1})-\util(t_{j+1})\bigr)\|_{\dot W(\ir)}\\
&\ \lesssim \|P_N e^{i(t-t_0)\Delta}\bigl(u(t_0)-\util(t_0)\bigr)\|_{\dot W(\ir)}+ \|P_N e\|_{\dot N^1(\ir)}\\
&\ \quad + \bigl\|P_N \bigl[(i\partial_t+\Delta)(u-\util)+e\bigr]\bigr\|_{\dot N^1(\ir)}.
\end{align*}
Squaring the above inequality, summing over all dyadic $N$'s, and using \eqref{closer},
\eqref{error small}, the dual of \eqref{square sum}, and the inductive hypothesis, we estimate
\begin{align*}
\Bigl(\sum_N&\|P_N  e^{i(t-t_{j+1})\Delta}\bigl(u(t_{j+1})-\util(t_{j+1})\bigr)\|_{\dot W(\ir)}^2\Bigr)^{1/2}\\
&\lesssim \Bigl(\sum_N\|P_N e^{i(t-t_0)\Delta}\bigl(u(t_0)-\util(t_0)\bigr)\|_{\dot W(\ir)}^2\Bigr)^{1/2}\\
&\quad +\Bigl(\sum_N\|P_N e\|_{\dot N^1(\ir)}^2\Bigr)^{1/2}\\
&\quad +\Bigl(\sum_N\bigl\|P_N \bigl[(i\partial_t+\Delta)(u-\util)+e\bigr]\bigr\|_{\dot N^1(\ir)}^2\Bigr)^{1/2}\\
&\lesssim \eps +\|e\|_{\dot N^1(\ir)}+\|(i\partial_t+\Delta)(u-\util)+e\|_{\dot N^1(\ir)} \\
&\lesssim \eps + \sum_{k=0}^j C(k)\bigl(\eps+\eps^{\frac{7}{(n-2)^2}}\bigr).
\end{align*}
Here, $C(k)$ depends only on $k$, $E$, $E'$, and $\eps_0$. Choosing $\eps_1$ sufficiently small depending on $N_1$, $E$,
and $E'$, we can continue the inductive argument.  Note that the final constants can easily be chosen to depend in a non-decreasing
manner on $E, E', M$ (which is quite plausible, given that increasing those parameters can only serve to worsen the situation).

This concludes the proof of Theorem~\ref{long-time theorem}. As a consequence of this theorem, we will derive scattering
results. Let us start by proving that a finite bound of the $L_{t,x}^{\frac{2(n+2)}{n-2}}$-norm of the solution to
\eqref{equation 1} implies scattering. Indeed, we have

\begin{corollary}[$L^p_{t,x}$ bounds imply scattering]\label{L^p implies scattering}
Let $u_0\in\dot{H}^1_x$ and let $u$ be a global solution to \eqref{equation 1} such that
\begin{align}\label{assume L^p}
\|u\|_{L_{t,x}^{\frac{2(n+2)}{n-2}}(\R\times\R^n)}\leq M
\end{align}
for some constant $M>0$. Then there exist finite energy solutions $u_{\pm}(t,x)$ to the free Schr\"odinger equation
$(i\partial_t+\Delta)u_{\pm}=0$ such that
$$
\|u_{\pm}(t)-u(t)\|_{\dot{H}^1_x}\rightarrow 0
$$
as $t\rightarrow \pm\infty$. Furthermore, the maps $u_0\mapsto u_{\pm}(0)$ are continuous from $\dot{H}^1_x$ to itself.
\end{corollary}

\begin{proof}
We will only prove the statement for $u_{+}$, since the proof for $u_{-}$ follows similarly.
Let us first construct the scattering state $u_{+}(0)$. For $t>0$ define $v(t) = e^{-it\Delta}u(t)$. We will show
that $v(t)$ converges in $\dot{H}^1_x$ as $t\rightarrow \infty$, and define $u_{+}(0)$ to be the limit.
Indeed, from Duhamel's formula \eqref{duhamel} we have
\begin{align}\label{v}
v(t) = u(0) - i\int_{0}^{t} e^{-is\Delta}f(u(s))ds.
\end{align}
Therefore, for $0<\tau<t$,
$$
v(t)-v(\tau)=-i\int_{\tau}^{t}e^{-is\Delta}f(u(s))ds.
$$
By \eqref{chain-game}, \eqref{lls-special-2}, and Lemma \ref{lemma strichartz norms}, we have
\begin{align*}
\|v(t)-v(\tau)\|_{\dot{H}^1_x}
&\lesssim \|f(u)\|_{\dot N^1([\tau,t]\times\R^n)}\\
&\lesssim \|u\|_{L_{t,x}^{\frac{2(n+2)}{n-2}}([\tau,t]\times\R^n)}^{\frac{4}{n-2}}\|u\|_{\dot W([\tau,t]\times\R^n)}\\
&\lesssim \|u\|_{L_{t,x}^{\frac{2(n+2)}{n-2}}([\tau,t]\times\R^n)}^{\frac{4}{n-2}}\|u\|_{\dot{S}^1([\tau,t]\times\R^n)}.
\end{align*}
However, \eqref{assume L^p} implies $\|u\|_{\dot{S}^1(\R\times\R^n)}\leq C(E,M)$ by the same argument as in
the proof of Theorem~\ref{long-time theorem}, where $E$ denotes the kinetic energy of the initial data $u_0$. Also by
\eqref{assume L^p}, for any $\eta>0$ there exists $t_{\eta}\in \R_{+}$ such that
$$
\|u\|_{L^{\frac{2(n+2)}{n-2}}_{t,x}([t,\infty)\times\R^n)}\leq \eta
$$
whenever $t>t_{\eta}$. Hence,
\begin{center}
$\|v(t)-v(\tau)\|_{\dot{H}^1_x}\rightarrow 0 \quad$ as $t,\tau\rightarrow \infty$.
\end{center}
In particular, this implies that $u_{+}(0)$ is well defined. Also, inspecting \eqref{v} one easily sees that
\begin{align}
u_{+}(0)=u_0- i\int_{0}^{\infty}e^{-is\Delta}f(u(s))ds
\end{align}
and thus
\begin{align}\label{u+}
u_{+}(t)=e^{it\Delta}u_0- i\int_{0}^{\infty}e^{i(t-s)\Delta}f(u(s))ds.
\end{align}
By the same arguments as above, \eqref{u+} and Duhamel's formula \eqref{duhamel} imply that
$\|u_{+}(t)-u(t)\|_{\dot{H}^1_x}\rightarrow 0$ as $t\rightarrow\infty$.

Similar estimates prove that the inverse wave operator $u_0\mapsto u_{+}(0)$ is continuous from $\dot{H}^1_x$ to
itself subject to the assumption \eqref{assume L^p} (in fact, we obtain a H\"older continuity estimate with this
assumption). We skip the details.
\end{proof}

\begin{remark}
If we assume $u_0\in H^1_x$ in Corollary~\ref{L^p implies scattering}, then similar arguments yield
scattering in $H^1_x$, i.e., there exist finite energy solutions $u_{\pm}(t,x)$ to
the free Schr\"odinger equation $(i\partial_t+\Delta)u_{\pm}=0$ such that
$$
\|u_{\pm}(t)-u(t)\|_{H^1_x}\rightarrow 0 \quad \text{as} \  t\rightarrow \pm\infty.
$$
\end{remark}

\begin{remark}
If we knew that the problem \eqref{equation 1} were globally wellposed for arbitrary $\dot{H}^1_x$ (respectively $H_x^1$)
initial data, then standard arguments would also give asymptotic completeness, i.e., the maps $u_0\mapsto u_{\pm}(0)$
would be homeomorphisms from $\dot{H}^1_x$ (respectively $H_x^1$) to itself. See for instance \cite{cazenave:book} for this
argument in the energy-subcritical case.
\end{remark}

As a consequence of Corollary~\ref{L^p implies scattering} and the global well-posedness theory for small initial data
(see Corollary~\ref{cor lwp}), we obtain scattering for solutions of \eqref{equation 1} with initial data small in the
energy-norm $\dot{H}^1_x$:

\begin{corollary}
Let $u_0\in H^1_x$ be such that
\begin{align*}
\|u_0\|_{\dot{H}^1_x}\lesssim \eta_0
\end{align*}
with $\eta_0$ as in Theorem~\ref{lwp} and let $u$ be the unique global solution to \eqref{equation 1}. Then there exist
finite energy solutions $u_{\pm}(t,x)$ to the free Schr\"odinger equation $(i\partial_t+\Delta)u_{\pm}=0$ such that
$$
\|u_{\pm}(t)-u(t)\|_{H^1_x}\rightarrow 0
$$
as $t\rightarrow \pm\infty$. Moreover, the maps $u_0\mapsto u_{\pm}(0)$ are continuous from $H_x^1$ to itself
(in fact, we have a H\"older continuity estimate).
\end{corollary}


\begin{thebibliography}{10}

\bibitem{borg:scatter}
J. Bourgain, \emph{Global well-posedness of defocusing 3D critical
NLS in the radial case}, JAMS \textbf{12} (1999), 145-171.

\bibitem{borg:book}
J. Bourgain, \emph{New global well-posedness results for nonlinear
Schr\"odinger equations}, AMS Publications, 1999.

\bibitem{cw0}
T. Cazenave, F.B. Weissler, \emph{Some remarks on the nonlinear Schr\"odinger equation in the critical case},
Nonlinear semigroups, Partial Differential Equations and Attractors, Lecture Notes in Math. \textbf{1394} (1989), 18--29.

\bibitem{cwI}
T. Cazenave, F.B. Weissler, \emph{Critical nonlinear Schr\"odinger
Equation}, Non. Anal. TMA \textbf{14} (1990), 807--836.

\bibitem{cazbook}
T. Cazenave, \emph{An introduction to nonlinear Schr\"odinger
equations}, Textos de M\'etodos Matem\'aticos \textbf{26},
Instituto de Matem\'atica UFRJ, 1996.

\bibitem{cazenave:book}
T. Cazenave, \textit{Semilinear Schr\"odinger equations,}
Courant Lecture Notes in Mathematics, 10.
American Mathematical Society, 2003.

\bibitem{gopher}
J. Colliander, M. Keel, G. Staffilani, H. Takaoka, T. Tao, \emph{Global well-posedness and
scattering in the energy space for the critical nonlinear
Schr\"odinger equation in $\R^3$}, preprint.

\bibitem{foschi}
D. Foschi, \emph{Inhomogeneous Strichartz estimates}, Journal of Hyperbolic Differential Equations, vol.2 no.1 (2005).

\bibitem{twounique}
G. Furioli, E. Terraneo, \emph{Besov spaces and unconditional well-posedness
for the nonlinear Schr\"odinger equation in $\dot{H}^s$}, Comm. in Contemp. Math.
\textbf{5}  (2003),  349--367.

\bibitem{FPT_NLSunique}
G. Furioli, F. Planchon, E. Terraneo, \emph{Unconditional well-posedness for semilinear Schr\"odinger equations in $H^s$},
Harmonic analysis at Mount Holyoke, (South Hadley, MA, 2001), 147--156.

\bibitem{gv:localreference}
J. Ginibre, G. Velo, \emph{ The global Cauchy problem for the nonlinear Schr\"odinger equation revisited},
Ann. Inst. H. Poincare' Anal. Non Line'aire \textbf{2} (1985),309--327.

\bibitem{glassey}
R.T. Glassey, \emph{On the blowing up of solutions to the Cauchy problem for nonlinear Schr\"odinger operators}, J. Math. Phys. \textbf{8} (1977), 1794--1797.

\bibitem{grillakis:scatter}
M. Grillakis, \emph{On nonlinear Schr\"odinger equations}, Comm.
Partial Differential Equations \textbf{25} (2000), no. 9-10,
1827--1844.

\bibitem{kato}
T. Kato, \emph{On nonlinear Schr\"odinger equations},  Ann. Inst. H. Poincare Phys. Theor.
\textbf{46}  (1987),  113--129.

\bibitem{katounique}
T. Kato, \emph{On nonlinear Schr\"odinger equations, II.  $H^s$-solutions and unconditional well-posedness}, J. d'Analyse. Math. \textbf{67}, (1995), 281--306.

\bibitem{tao:keel}
M. Keel, T. Tao, \emph{Endpoint Strichartz Estimates}, Amer. Math.
J. \textbf{120} (1998), 955--980.

\bibitem{keraani}
S. Keraani, \emph{On the defect of compactness for the Strichartz estimates of the Schr\"odinger equations}, J. Diff. Eq. \textbf{175}, (2001), 353--392.

\bibitem{nakanishi}
K. Nakanishi, \emph{Scattering theory for nonlinear Klein-GOrdon equation with Sobolev critical power}, IMRN \textbf{1} (1999), 31--60.

\bibitem{rv}
E. Ryckman, M. Visan, \emph{Global well-posedness and scattering for the defocusing energy-critical nonlinear Schr\"odinger equation in $\R^{1+4}$}, preprint.

\bibitem{tao:gwp radial}
T. Tao, \emph{Global well-posedness and scattering for the
higher-dimensional energy-critical non-linear Schr\"odinger
equation for radial data}, to appear New York Journal of Math.

\bibitem{monica-thesis}
M. Visan, \emph{The defocusing energy-critical nonlinear Schr\"odinger equation in higher dimensions}, in preparation.

\end{thebibliography}
\end{document}